\documentclass[10pt,a4paper,titlepage]{article}
\usepackage[latin1]{inputenc}
\usepackage{amsmath}
\usepackage[ps,all,arc,rotate]{xy}
\usepackage{amsfonts}
\usepackage{amssymb}
\usepackage{amsthm}
\usepackage[dvips]{graphicx}
\usepackage{fancyhdr}
\usepackage{float}
\usepackage{young}
\author{Dawid Kielak}
\title{Low-dimensional free and linear representations of $\Out{3}$}
\makeatletter
\def\imod#1{\allowbreak\mkern10mu({\operator@font mod}\,\,#1)}
\makeatother
\numberwithin{figure}{section}
\theoremstyle{plain}
\newtheorem{thm}{Theorem}[section]
\newtheorem*{prop*}{Proposition}
\newtheorem*{thm*}{Theorem}
\newtheorem*{thmA}{Theorem \ref{result A}}

\newtheorem*{thmC}{Theorem \ref{result C}}
\newtheorem{prop}[thm]{Proposition}
\newtheorem{lem}[thm]{Lemma}

\theoremstyle{definition}
\newtheorem{ex}[thm]{Example}

\newtheorem{dfn}[thm]{Definition}
\newtheorem*{dfn*}{Definition}
\newtheorem{notation}[thm]{Notation}
\theoremstyle{remark}
\newtheorem{rmk}[thm]{Remark}

\def\C{\mathbb{C}}

\def\Z{\mathbb{Z}}

\def\K{\mathbb{K}}
\def\Q{\mathbb{Q}}

\def\Wlog{Without loss of generality }

\def\rep{representation }
\def\reps{representations }
\def\iff{if and only if }

\def\GL{\mathrm{GL}}
\def\SL{\mathrm{SL}}
\def\I{\mathrm{I}}

\def\into{\hookrightarrow}
\def\Torelli{\overline{\mathrm{IA}}}

\newcommand{\Out}[1]{\mathrm{Out}(F_{#1})}
\newcommand{\Aut}[1]{\mathrm{Aut}(F_{#1})}

\newdir{ >}{{}*!/-5pt/@{>}}
\newdir{< }{{}*!/+5pt/@{<}}
\def\s-{\smallsetminus}
\begin{document}
\textsc{\begin{LARGE}\begin{center} Low-dimensional free and linear representations of $\Out{3}$
\end{center}\end{LARGE}}

\medskip

\begin{center}
Dawid Kielak\footnotemark

University of Bonn

\texttt{kielak@math.uni-bonn.de}

\today
\end{center}

\medskip

\begin{center}
\begin{minipage}{0.65\textwidth}
\textsc{Abstract.} We study homomorphisms from $\Out{3}$ to $\Out{5}$, and $\GL_m(\K)$ for $m \leqslant 6$, where $\K$ is a field of characteristic other than 2 or 3. We conclude that all $\K$-linear representations of dimension at most 6 of $\Out{3}$ factor through $\GL_3(\Z)$, and that all homomorphisms $\Out{3} \to \Out{5}$ have finite image.
\end{minipage}
\end{center}

\bigskip

\footnotetext{The author was supported by the EPSRC of the United Kingdom}

\section{Introduction}

This paper constitutes a part of a project of understanding homomorphisms $\Out{n} \to \GL_m(\K)$ and $\Out{n} \to \Out{m}$, i.e. the $\K$-linear representation theory and the \emph{free representation theory} of $\Out{n}$ respectively (see~\cite{kielak2013}). The first problem, as well as linear representations of $\Aut{n}$,  has been studied for example by Potapchik--Rapinchuk~\cite{potapchikrapinchuk2000} and
Grunewald--Lubotzky~\cite{grunewaldlubotzky2006};
the latter has been addressed for example by Khramtsov~\cite{khramtsov1990}, Bogopolski--Puga~\cite{bogopol'skiipuga2002}, Bridson--Vogtmann~\cite{bridsonvogtmann2003, bridsonvogtmann2011}, and Aramayona--Leininger--Souto~\cite{aramayonaetal2009}.

The results obtained so far have a tendency of working only for large $n$; in particular the case of $\Out{3}$ has not been studied extensively.
It was known (thanks to a result of Khramtsov~\cite{khramtsov1990}) that there are no embeddings \[\Out{3} \into \Out{4} \] It was also known that there is an embedding \[\Out{3} \into \Out{55}\] (this is due to Bogopolski--Puga~\cite{bogopol'skiipuga2002}).

Our attempt to expand the understanding of free representations of $\Out{3}$ follows the same general route as our proof in~\cite{kielak2013}. First we investigate the linear representation theory of $\Out{3}$, and we prove
\begin{thmA}
Let $\phi$ be a $\K$-linear, six-dimensional \rep of $\Out{3}$, where $\K$ is a field of characteristic other than 2 or 3. Then $\phi$ factors as
\[ \xymatrix{ \Out{3} \ar[r]^\phi \ar[dr]^\pi & \GL_6(\K) \\ & \GL_3(\Z) \ar[u]}  \]
where $\pi$ is the natural projection $\Out{3} \to \mathrm{Out}(H_1(F_3,\Z)) \cong \GL_3(\Z)$.
\end{thmA}

\hyphenation{di-men-sion-al}

It is worth noting here that Turchin--Willwacher~\cite{TurchinWillwacher2015} constructed a 7-dimensional $\Q$-linear representation of $\Out{3}$ which does not factor through $\pi$.
Thus our result completes the search for the smallest dimension in which such a representation occurs (at least for fields of characteristic zero).

The second part of this paper deals with the action of two finite subgroups of $\Out{3}$ on graphs of rank 5. This was also our strategy in~\cite{kielak2013}, yet the situation here is quite different. The finite groups under consideration do not contain (comparatively) large simple groups, as is the case in the higher rank case; on the other hand, the groups are of order 48, and hence are rather tangible. Our considerations yield
\begin{thmC}
Suppose $\phi : \Out{3} \to \Out{5}$ is a homomorphism. Then the image of $\phi$ is finite.
\end{thmC}
The general question of finding the smallest $m>3$ such that there is a homomorphism $\Out{3} \to \Out{m}$ that is injective, or at least has infinite image, remains open.

\noindent \textbf{Acknowledgements.}
The author wishes to thank Martin R. Bridson for all his help.

\section{Notation and preliminaries}
\label{secnotation}
Let us first establish some conventions and definitions.

\begin{dfn}
\label{graphdef}
We say that $X$ is a \emph{graph} \iff it is a 1-dimensional CW complex. The closed 1-cells of $X$ will be called \emph{edges}, the 0-cells will be called \emph{vertices}. The sets of vertices and edges of a graph will be denoted by $V(X)$ and $E(X)$ respectively. The points of intersection of an edge with the vertex set are referred to as \emph{endpoints} of the edge.

We will equip $X$ with the standard path metric in which the length of each edge is 1.

Given two graphs $X$ and $Y$, a function $f : X \to Y$ is a \emph{morphism of graphs}  \iff $f$ is a continuous map sending $V(X)$ to $V(Y)$, and sending each open edge in $X$ either to a vertex in $Y$ or isometrically onto an open edge in $Y$. Note that a morphism can invert edges.

When we say that a group $G$ \emph{acts on a graph} $X$, we mean that it acts by graph morphisms.

We say that a graph $X$ is \emph{directed} \iff it comes equipped with a map $o : E(X) \to X$ such that $o(e)$ is a point on the interior of $e$ of distance $\frac13$ from one of its endpoints. We also define $\iota, \tau : E(X) \to V(X)$ by setting $\tau(e)$ to be the endpoint of $e$ closest to $o(e)$, and $\iota(e)$ to be the endpoint of $e$ farthest from $o(e)$. Note that we allow $\iota(e) = \tau(e)$.

%
The \emph{rank} of a connected graph is defined to be the size of a minimal generating set of its fundamental group (which is a free group). 
\end{dfn}

Let us also define two families of graphs.

\begin{dfn}
The graph with one vertex and $n$ edges will be referred to as the \emph{$n$-rose}.

The graph with two vertices and $n$ edges, such that each edge has two distinct endpoints, will be referred to as the \emph{$n$-cage}.
\end{dfn}

%
%
%
\begin{notation}
Let $G$ be a group. We will adopt the following notation:
\begin{itemize}
\item for two elements $g,h \in G$, we define $g^h=h^{-1}gh$;
\item for two elements $g,h \in G$, we define $[g,h]=ghg^{-1}h^{-1}$;
\end{itemize}
We will also use $\Z_k$ to denote the cyclic group of order $k$.
\end{notation}

\begin{dfn}
Let us introduce the following notation for elements of $\mathrm{Aut}(F_n)$, the automorphism group of $F_n$, where $F_n$ is the free group on $\{ a_1, \ldots, a_n\}$:
\begin{eqnarray*}
\epsilon_i : \left\{ \begin{array}{cccl} a_i & \mapsto & a_i^{-1}, \\ a_j & \mapsto & a_j, & j \neq i \end{array} \right. & &
\sigma_{ij} : \left\{ \begin{array}{cccl} a_i & \mapsto & a_j, \\ a_j & \mapsto & a_i, & \\ a_k & \mapsto & a_k, & k \not \in \{i,j\} \end{array} \right. \\
\rho_{ij} : \left\{ \begin{array}{cccl} a_i & \mapsto & a_i a_j, \\ a_k & \mapsto & a_k, &  k \neq i \end{array} \right. & &
\lambda_{ij} : \left\{ \begin{array}{cccl} a_i & \mapsto & a_j a_i, \\ a_k & \mapsto & a_k, &  k \neq i \end{array} \right.
\end{eqnarray*}
Let us also define $\Delta = \prod_{i=1}^n \epsilon_i$ and
\begin{displaymath}
\sigma_{i (n+1)} : \left\{ \begin{array}{cccl} a_i & \mapsto & a_i^{-1}, \\ a_j & \mapsto & a_j a_i^{-1}, & j \neq i \end{array} \right.
\end{displaymath}
\end{dfn}

Below we give an explicit presentation of $\Out{n}$, the outer automorphism group of $F_n$:

\begin{thm}[Gersten's presentation~\cite{gersten1984}]
\label{gerstenpres}
Suppose $n \geqslant 3$. The group $\Out{n}$ is generated by $\{ \epsilon_1, \rho_{ij}, \lambda_{ij} \mid i,j = 1,\ldots,n , \, i \neq j \}$, with relations
\begin{itemize}
\item $[\rho_{ij},\rho_{kl}]=[\lambda_{ij},\lambda_{kl}]=1$ for $k \not \in \{i,j\}, l \neq i$;
\item $[\lambda_{ij},\rho_{kl}]=1$ for $k \neq j, l \neq i$;
\item $[\rho_{ij}^{-1},\rho_{jk}^{-1}]=[\rho_{ij},\lambda_{jk}]=[\rho_{ij}^{-1},\rho_{jk}]^{-1}=[\rho_{ij},\lambda_{jk}^{-1}]^{-1}=\rho_{ik}^{-1}$ for $k \not \in \{i,j\}$;
\item $[\lambda_{ij}^{-1},\lambda_{jk}^{-1}]=[\lambda_{ij},\rho_{jk}]=[\lambda_{ij}^{-1},\lambda_{jk}]^{-1}=[\lambda_{ij},\rho_{jk}^{-1}]^{-1}=\lambda_{ik}^{-1}$ for $k \not \in \{i,j\}$;
\item $\rho_{ij} \rho_{ji}^{-1} \lambda_{ij}=\lambda_{ij} \lambda_{ji}^{-1} \rho_{ij}, \ (\rho_{ij} \rho_{ji}^{-1} \lambda_{ij})^4=1$;
\item $[\epsilon_1, \rho_{ij}] = [\epsilon_1, \lambda_{ij}]=1$ for $i,j \neq 1$;
\item $\rho_{12}^{\epsilon_1}=\lambda_{12}^{-1}, \ \rho_{21}^{\epsilon_1}=\rho_{21}^{-1}$;
\item $\epsilon_1^2=1$;
\item $\prod_{i\neq j} \rho_{ij} \lambda_{ij}^{-1} = 1$ for each fixed $j$.
\end{itemize}
\end{thm}

Note the action of $\Aut{n}$ on $F_n$ and $\Out{n}$ on the conjugacy classes of $F_n$ is \textbf{on the left}.

\begin{dfn}
Let us define some finite subgroups of $\Out{n}$:
\begin{eqnarray*}
S_n &\cong&\langle \{\sigma_{ij} \mid i,j = 1,\ldots,n, \ i \neq j\} \rangle\\
S_{n+1} &\cong&\langle \{\sigma_{ij} \mid i,j = 1,\ldots,n+1, \ i \neq j \} \rangle\\
\Z_2^n \rtimes S_n\cong W_n& =& \langle \{ \epsilon_1, \sigma_{ij} \mid i,j=1\ldots,n, \ i \neq j \} \rangle \\
\Z_2 \times S_{n+1}\cong G_n & =& \langle \{ \Delta, \sigma_{ij} \mid i,j=1\ldots,n+1, \ i\neq j \} \rangle
\end{eqnarray*}
\end{dfn}

We do not give distinctive names to the first two (symmetric) groups; instead, we will usually refer to them as respectively $S_n < W_n$ and $S_{n+1}<G_n$. More generally, whenever we mention $S_n$ or $S_{n+1}$ as subgroups of $\Out{n}$, we mean these two groups.

In the case of $\Out{3}$, define $V_4$ and $A_4$ to be the Klein 4-group and the alternating group of degree 4 satisfying
\[ V_4 < A_4 < S_4 < G_3 < \Out{3}. \]

Note that we abuse notation by also using $S_n$ to denote the abstract symmetric group of degree $n$, and $A_n < S_n$ to denote its maximal alternating subgroup.

Note that, if $i,j \leqslant n$, we have
\begin{displaymath}
\epsilon_i \sigma_{ij} = \lambda_{ij} \lambda_{ji}^{-1} \rho_{ij}=\rho_{ij} \rho_{ji}^{-1} \lambda_{ij}
\end{displaymath}
and
\[ \epsilon_1 \sigma_{i (n+1)} = \prod_{j \neq i} \rho_{ji} = \prod_{j \neq i} \lambda_{ji} \]
and the subgroup ${S_n < \Out{n}}$ defined above acts on the sets
\begin{align*}
\{ \epsilon_i & \mid i=1\ldots,n\}, \\ \{ \rho_{ij} & \mid i,j=1\ldots,n, \ i \neq j\}, \textrm{ and} \\ \{ \lambda_{ij} & \mid i,j=1\ldots,n, \ i \neq j\}
\end{align*}
by permuting the indices in the natural way.

%
%
\section{Representations of $\Out{3}$ in dimensions five and six}

Before we start investigating 5- and 6-dimensional representations of $\Out{3}$, let us first prove the following.

\begin{prop}
\label{killV_4}
Suppose $\phi: \Out{3} \to G$ is a group homomorphism such that its kernel contains $V_4$. Then $\phi$ factors as
\[ \phi: \Out{3} \to \Z_2 \to G\]
and the map is determined by the image of $\epsilon_1$.
\begin{proof}
Since $V_4$ lies in the kernel, we have $\phi(\sigma_{14}) = \phi(\sigma_{23})$. Hence
\[\phi(\rho_{21}) = \phi(\rho_{21}^{\rho_{21} \rho_{31}}) = \phi(\rho_{21}^{\epsilon_1 \sigma_{14}}) = \phi(\rho_{21})^{ \phi(\epsilon_1) \phi( \sigma_{23}) } = \phi(\rho_{31}^{-1})\]
and so
\[ \phi(\rho_{21}) = \phi(\rho_{21}^{\epsilon_3}) = \phi(\rho_{31}^{-1})^{\phi(\epsilon_3)} = \phi(\lambda_{31})  \]
Now
\[ \phi(\rho_{31}^{-1}) = \phi([\rho_{32}^{-1},\rho_{21}^{-1}]) = [\phi(\rho_{32}^{-1}),\phi(\lambda_{31}^{-1})]=1 \]
Thus $\rho_{31}$ lies in the kernel of $\phi$. We can however conjugate $\rho_{31}$ to each $\rho_{ij}$ using $S_3$, and so all elements $\rho_{ij}$ lie in the kernel. The result follows.
\end{proof}
\end{prop}

Let us recall some basic terminology of the representation theory of symmetric groups.
\begin{rmk}
Let $S_n$ be a symmetric group of rank $n$, and let $\K$ be a field with $\mathrm{char}(\K) = 0$ or $\mathrm{char}(\K) > n$. Then there is a one-to-one correspondence between irreducible $\K$-linear \reps of $S_n$ and partitions of $n$. In particular the \rep corresponding to
\begin{itemize}
\item $(n)$ is called \emph{trivial};
\item $(1^n)$ is called \emph{determinant};
\item $(n-1,1)$ is called \emph{standard};
\item $(2,1^{n-1})$ is called \emph{signed standard}.
\end{itemize}
Moreover, a direct sum of the standard and trivial \rep is known as the \emph{permutation} representation, and a direct sum of the signed standard and determinant \rep is known as the \emph{signed permutation} representation.
\end{rmk}

Now let us turn our attention to representations of the group $W_n$.

\begin{dfn}
\label{repofW_n}
Let $V$ be a representation of $W_n$. Let $N = \{1,\ldots,n\}$. Define
\begin{itemize}
\item for each $I \subseteq N$,  $E_I = \{ v \in V \mid \epsilon_i v=(-1)^{\chi_I(i)}v\}$, where $\chi_I$ is the characteristic function of $I$;
\item $V_i = \bigoplus_{\vert I \vert =i} E_I$
\end{itemize}
\end{dfn}

We will slightly abuse notation, and often omit parentheses and write $E_1$ for $E_{\{1\}}$, etc.

\begin{lem}
\label{dimensionofV_i}
Let $V$ be a representation of $W_n$. Then $\dim V_i = {n \choose i} \dim E_I $ where $\vert I \vert = i$.
\begin{proof}
The symmetric group $S_n < W_n$ acts on $\{ \epsilon_1, \ldots , \epsilon_n \}$ by permuting the indices in the natural way. Hence its action on $V_i$ will permute subspaces $E_I$ by permuting subsets of $N$ of size $i$. Thus each $E_I$, for a fixed size of $I$, has the same dimension. The result follows.
\end{proof}
\end{lem}

An immediate consequence of the above is the following.
\begin{lem}
\label{repsofdim2}
\label{lem: 2-dim reps}
Let $V$ be a 2-dimensional $\K$-linear representation of $\Out{3}$, where $\mathrm{char}(\K) \neq 2$. Then the representation factors as \[\Out{3} \to \Z_2 \to \GL_2(\K)\] and is determined by the image of $\epsilon_1$.
\begin{proof}
Let $\phi: \Out{3} \to \GL_2(\K)$ be the representation. Lemma~\ref{dimensionofV_i} tells us that (with the notation of Definition~\ref{repofW_n}) $V = V_0 \oplus V_3$.

Since $S_4 < G_3$ commutes with $\Delta$, $V = V_0 \oplus V_3$ is a decomposition of $S_4$-modules. Now each of these submodules has dimension at most 2. There are at most three irreducible $\K$-linear representations of $S_4$ of dimension at most 2: the trivial representation (corresponding to partition $(4)$), the determinant representation (corresponding to partition $(1^4)$), and the one given by a partition $(2,2)$  (note that the latter might not be irreducible when $\mathrm{char}(\K) =3$).
In all three cases, the action of $V_4 < S_4$ is trivial. This implies that we have satisfied all the requirements of Proposition~\ref{killV_4}, and the result follows.
\end{proof}
\end{lem}

\begin{lem}
\label{diamonds}
Let $V$ be a $\K$-linear \rep of $\Out{n}$, where $\K$ is a field of characteristic other than 2. Then, with the notation above, we have
\[ V = \bigoplus_{i=0}^{n} V_i\]
and for each $i \neq j$, $J \subseteq N \smallsetminus \{i,j\}$ we have
\[ \rho_{ij} (E_J \oplus E_{J \cup \{i\}} \oplus E_{J \cup \{j\}} \oplus E_{J \cup \{i,j\}}) = E_J \oplus E_{J \cup \{i\}} \oplus E_{J \cup \{j\}} \oplus E_{J \cup \{i,j\}}.\]
An identical statement holds for $\lambda_{ij}$.
\begin{proof}
The first statement follows directly from the fact that we can simultaneously diagonalise commuting involutions $\epsilon_i$, since we are working over a field $\K$ whose characteristic is not 2.

For the second statement, let us note that $ [\rho_{ij},\epsilon_k] = 1$ for each $k \not\in \{i,j\}$. Hence for each $I \subseteq N$:
\[ \rho_{ij} (E_I) \leqslant \bigoplus_{J \triangle I \subseteq \{i,j\}} E_J \]
where $A \triangle B$ denotes the symmetric difference of two sets $A$ and $B$. An identical argument works for $\lambda_{ij}$.
\end{proof}
\end{lem}

To help us visualise the combinatorics of \reps of $\Out{n}$ we are going to use the following diagrams.

\begin{dfn}
Suppose $V$ is a finite dimensional, $\K$-linear representation of $\Out{n}$ over any field $\K$, and let $x \in \Out{n}$. Let us use the notation of Definition~\ref{repofW_n}. We define \emph{the minimal diagram for $x$ over $V$} (often abbreviated to the minimal diagram for $x$) to be a directed graph $D$ with the vertex set equal to a subset $S$ of the power set of $N = \{1 ,2 ,\ldots, n\}$, where $I \in S$ \iff $E_I \neq \{0\}$, and the edge set given by the following rule: there is a directed edge from $I$ to $J$ \iff $p_J\big(x(E_I)\big) \neq \{ 0 \}$, where $p_J : \bigoplus_{K \subseteq N} E_K \to E_J$ is the natural projection.

We also say that a graph $D'$ is \emph{a diagram for $x$ over $V$} \iff the minimal diagram $D$ for $x$ is a subgraph of $D'$.
\end{dfn}

In practice, when realising these diagrams in terms of actual pictures, we are going to align vertices corresponding to subsets of $N$ of the same cardinality in horizontal lines; each such line will correspond to some $V_i$. We are also going to represent edges as follows: if two vertices are joined by two directed edges, we are going to draw one edge without any arrowheads between them; we are not going to draw edges from a vertex to itself -- instead, if a vertex does not have such a loop, then all edges emanating from it will be drawn with a tail (see example below);
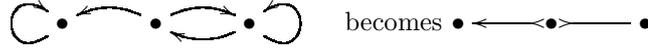
\[ \xymatrix{ \bullet \ar@(dl,ul)[] & \bullet \ar@/_/[l] \ar@/^/[r] & \bullet \ar@/^/[l] \ar@(dr,ur)[] & \textrm{ becomes }  \bullet & \bullet \ar@{ >->}[l] \ar@{ >-}[r] & \bullet
} \]

To get a firmer grip on these diagrams, let us have a look at a number of facts one can easily deduce from (not necessarily minimal) diagrams.

\begin{ex}
\label{ex: component}
Let $\Gamma_0$ be a connected component of $\Gamma$, a diagram for $x$. Let $ v \in \bigoplus_{I \in V(\Gamma_0)} E_I$ be a vector. Then $v = \sum v_I$ where $v_I \in E_I$.
Let $J \not\in V(\Gamma_0)$. Note that there are no edges between $J$ and $V(\Gamma_0)$, and so $p_J \big(x(v_I)\big) = 0$ for all $I \in V(\Gamma_0)$. Hence $x(v_I) \in \bigoplus_{I \in V(\Gamma_0)} E_I$ and therefore
\[ x(\bigoplus_{I \in V(\Gamma_0)} \!\!\!\!\! E_I) = \bigoplus_{I \in V(\Gamma_0)} \!\!\!\!\! E_I \]
\end{ex}

The following illustrates the relationship between our diagrams and matrices.

\begin{ex}
\label{ex: matrices}
Suppose we have a diagram for $x$ with two vertices, $I$ and $J$ say, such that the union of the connected components containing these vertices does not contain any other vertex. Fix a basis for $E_I$ and $E_J$. The following illustrates
the way the $x$ action on $E_I \oplus E_J$ (seen as a matrix) depends on the diagram:
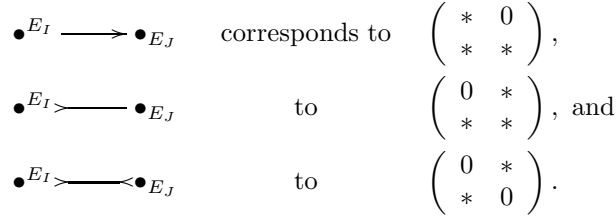
\begin{eqnarray*}
\xymatrix{ \bullet^{E_I} \ar[r] & \bullet_{E_J}} & \textrm{ corresponds to } & \left( \begin{array}{cc} * & 0 \\ * & * \end{array} \right), \\
\xymatrix{ \bullet^{E_I} \ar@{ >-}[r] & \bullet_{E_J}} & \textrm{ to } & \left( \begin{array}{cc} 0 & * \\ * & * \end{array} \right), \textrm{ and} \\
\xymatrix{ \bullet^{E_I} \ar@{ >-< }[r] & \bullet_{E_J}} & \textrm{ to } & \left( \begin{array}{cc} 0 & * \\ * & 0 \end{array} \right).
\end{eqnarray*}
\end{ex}

\begin{ex}
\label{ex: away arrow}
Suppose we have a diagram $\Gamma$ for $x$ such that $\Gamma$ has a connected component with only two vertices, $I$ and $J$ say, as depicted below.
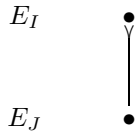
\[ \xymatrix{ E_I & \bullet \ar@{ >-}[d] \\ E_J & \bullet }\]
Exercise~\ref{ex: component} tells us that $E_I \oplus E_J$ is $x$-invariant, and thus $x\vert_{E_I \oplus E_J}$ is an isomorphism. Let $\{ v_1, \ldots, v_k \}$ be a basis for $E_I$. Our diagram tells us that $p_I\big(x(v_i)\big)=0$ for each $i$, and so $x(v_i) \in E_J$. Since $x$ is an isomorphism, we immediately see that \[ \{ x(v_1),\ldots,x(v_k)\} \] is a linearly independent set, and hence $\dim E_J \geqslant \dim E_I$.
\end{ex}

\begin{ex}
\label{ex: W_n action}
Let $\Gamma$ be a diagram for $x$. Note that $p_J\big(x^{\epsilon_i}(v)\big)=\pm p_J \big( x(v) \big)$ whenever $v \in E_I$ for some $I$. Therefore $p_J\big(x^{\epsilon_i}(E_I)\big)=\{ 0 \}$ \iff $p_J\big(x(E_I)\big)=\{ 0 \}$, and so $\Gamma$ is also a diagram for $x^{\epsilon_i}$.

Now consider $\sigma \in S_n$. We have
\[ p_J\big(x^{\sigma}(E_I)\big)= \sigma^{-1} \Big(p_{\sigma(J)}\big(x(E_{\sigma(I)})\big)\Big)\]
and so the image of $\Gamma$ under the graph morphism induced by $\I \mapsto \sigma(I)$ is a diagram for $x^\sigma$.
\end{ex}

We will use the last example very often, for example to relate diagrams for $\rho_{21}$ with ones for $\rho_{21}^{-1} = \rho_{21}^{\epsilon_1}$ or $\rho_{31} = \rho_{21}^{\sigma_{23}}$.

\begin{lem}
\label{lem: annoying case}
Let $\phi : \Out 3 \to V$ be a representation such that
we have a diagram for $\rho_{21}$ of the form
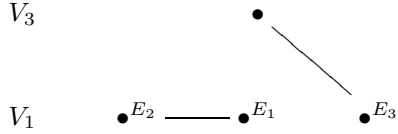
\[ \xymatrix{ V_3 & & \bullet \ar@{-}[dr] \\ V_1 & \bullet^{E_2} \ar@{-}[r] & \bullet^{E_1} & \bullet^{E_3} } \]
where $\dim E_i = 1$ for all $i$. Then $\rho_{21}$ has a diagram
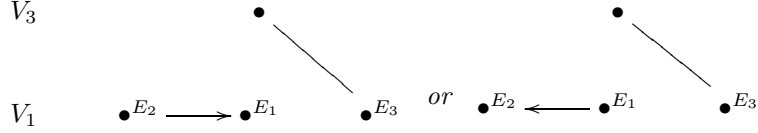
\[ \xymatrix{ V_3 & & \bullet \ar@{-}[dr] \\ V_1 & \bullet^{E_2} \ar[r] & \bullet^{E_1} & \bullet^{E_3}}
\xymatrix{ \\ \textrm{ or }}
\xymatrix{  & \bullet \ar@{-}[dr] \\ \bullet^{E_2}  & \bullet^{E_1} \ar[l] & \bullet^{E_3}  } \]
\end{lem}

Note that the above diagrams are not necessarily minimal; this means that the lines in the top diagram between $E_1$ and $E_2$, and between $V_3$ and $E_3$, may represent the real situation, or it may be possible to replace each by an edge directed from the right to the left, or one directed from the left
to the right, or no edge at all. Lemma 3.12 says that the line between $E_2$ and $E_1$ can indeed be replaced by either one of the directed edges or no edge at all.

\begin{proof}

Firstly, note that the left-hand side diagram corresponds to
\[{p_2 \big( \rho_{21} (E_1) \big) = \{ 0 \} }\]
 and the right-hand side diagram to $p_1 \big( \rho_{21} (E_2) \big) = \{ 0 \}$.
Suppose for a contradiction that we have neither of the diagrams, that is that $p_1 \big( \rho_{21} (E_2) \big) \neq \{ 0 \}$ and ${p_2 \big( \rho_{21} (E_1) \big) \neq \{ 0 \} }$. We claim that then $\rho_{21}$ has a diagram
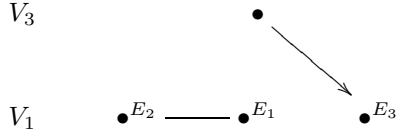
\[ \xymatrix{ V_3 & & \bullet \ar[dr] \\ V_1 & \bullet^{E_2} \ar@{-}[r] & \bullet^{E_1} & \bullet^{E_3} } \]

\smallskip

Once we have proven the above claim, take $x \in E_2 \s- \{0\}$. Then, by the assumptions of our claim, $p_1 \big( \rho_{21}(x) \big) \neq 0$, since $E_2$ is 1 dimensional, and so spanned by $x$.
Thus we have
\[ \rho_{21}(x) = x_1 + x_2 \]
with $x_1 \in E_1 \s- \{0\}$ and $x_2 \in E_2$.
Now
\[ p_3 \big( \rho_{31} \rho_{21}(x) \big) = p_3 \big( \rho_{31} (x_1 ) \big) + p_3 \big( \rho_{31} (x_2) \big) \]
By the conclusion of our claim, the following is a diagram for $\rho_{21}$
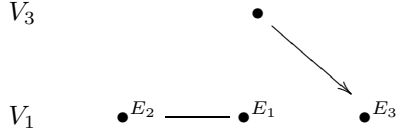
\[ \xymatrix{ V_3 & & \bullet \ar[dr] \\ V_1 & \bullet^{E_2} \ar@{-}[r] & \bullet^{E_1} & \bullet^{E_3} } \]
By Exercise~\ref{ex: W_n action} (using the action of $\sigma_{23}$) the following is then a diagram for $\rho_{31}$
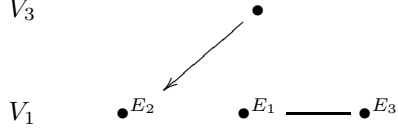
\[ \xymatrix{ V_3 & & \bullet \ar[dl] \\ V_1 & \bullet^{E_2}  & \bullet^{E_1} \ar@{-}[r] & \bullet^{E_3} } \]
Therefore
\[ p_3 \big( \rho_{31} (x_2) \big) = 0 \]
and so
\[ p_3 \big( \rho_{31} \rho_{21}(x) \big) = p_3 \big( \rho_{31} (x_1 ) \big) \]

Note that we can apply the assumptions of our claim to $\rho_{31}$ as well, using the relation $\rho_{31} = {\rho_{21}}^{\sigma_{23} }$ and Example~\ref{ex: W_n action}; specifically we may assume that $p_3 \big( \rho_{31}(E_1) \big) \neq \{0\}$.
Now, observing that
$x_1$ spans $E_1$ (which is 1 dimensional), we conclude that
\[ p_3 \big( \rho_{31} \rho_{21}(x) \big) \neq 0 \]

But $\rho_{31}(x) \in E_2$ by the diagram above (recall that $x \in E_2$), and thus $\rho_{21} \rho_{31} (x) \in E_2 \oplus E_1$, which in turn implies that
\[ p_3 \big( \rho_{21} \rho_{31}(x) \big) = 0 \]
This contradicts the relation $[\rho_{21},\rho_{31}]=1$, and our proof is complete.

\smallskip

Now, to prove the claim,
let $v_1 \in E_1 \smallsetminus \{ 0 \}$. Then $v_3 = \sigma_{23} p_2 \big( \rho_{21} (v_1) \big)$ generates $E_3$. Now if \[p_{1,2,3} \big( \rho_{21}(v_3) \big) = 0\] then we have proven our claim. If not, let $U = \langle u \rangle$ be a subspace of $V_3$ of dimension 1, where $u=p_{1,2,3} \big( \rho_{21}(v_3) \big)$.
Note that $\rho_{21}v_3 \in u + E_3$, and so
\[ \rho_{21}^{-1} (-v_3) = \epsilon_1 \rho_{21} \epsilon_1  (-v_3) \in u + E_3 .\]
In particular, $\rho_{21}^{-1} (-v_3) - v_3' = u$ for some $v_3' \in E_3$.
Now
\[\rho_{21} u = \rho_{21} \big( \rho_{21}^{-1} (-v_3) - v_3' \big) \in E_3 \oplus U , \] and hence $U \oplus E_3$ is $\rho_{21}$-invariant.

Let us rewrite $\rho_{21} \rho_{31} = \rho_{31} \rho_{21}$ as
\[ \rho_{21} \sigma_{23} \rho_{21} \sigma_{23} = \sigma_{23} \rho_{21} \sigma_{23} \rho_{21}\]
which yields $[\rho_{21} \sigma_{23}\rho_{21}, \sigma_{23}]=1$.

Now $E_1$ is $\sigma_{23}$-invariant and one dimensional, and so each non-zero vector in $E_1$ is an eigenvector of $\sigma_{23}$ with eigenvalue $\mu = \pm 1$. In particular $v_1$ lies in the $\mu$-eigenspace of $\sigma_{23}$, and hence so does $\rho_{21} \sigma_{23} \rho_{21} (v_1)$, since $\sigma_{23}$ and $\rho_{21} \sigma_{23}\rho_{21}$ commute.
Since $V_3$ and $V_1$ are $\sigma_{23}$-invariant (in fact $S_3$-invariant), we also have $p_{1,2,3} \big( \rho_{21} \sigma_{23} \rho_{21} (v_1) \big)$ lying in the $\mu$-eigenspace of $\sigma_{23}$.

Now let
\[ \rho_{21}(v_1) = y_1 + y_2 \]
with $y_i \in E_i$. We have
\begin{eqnarray*} p_{1,2,3} \big( \rho_{21} \sigma_{23} \rho_{21} (v_1) \big) &=& p_{1,2,3} \big( \rho_{21} \sigma_{23} (y_1+y_2) \big) \\ &=& p_{1,2,3} \big( \rho_{21} (\mu y_1 + \sigma_{23}  (y_2)) \big) \\ &=& p_{1,2,3} \big( \rho_{21} \sigma_{23} (y_2) \big)\end{eqnarray*}
since $\rho_{21}(\mu y_2) \in E_1 \oplus E_2$ by assumption of our lemma, and thus lies in the kernel of $p_{1,2,3}$.

By definition we have
\[ v_3 = \sigma_{23} p_2 \rho_{21} (v_1) = \sigma_{23} (y_2)\]
and so we deduce
\[ p_{1,2,3} \big( \rho_{21} \sigma_{23} \rho_{21} (v_1) \big) = p_{1,2,3} \big( \rho_{21} (v_3) \big) = u\]
where the last equality is the definition of $u$.
 Therefore $U$ lies in the $\mu$-eigenspace of $\sigma_{23}$.

Note that the eigenspaces of $\Delta$ are $V_0 \oplus V_2$ and $V_1 \oplus V_3$, and on each $\Delta$ acts as $\pm 1$. Hence, since $V_0 \oplus V_2 = \{0\}$, we see that $[\phi(\rho_{21}),\phi(\Delta)]=1$. Therefore, since $\rho_{ij}^\Delta = \lambda_{ij}$, the elements $\rho_{ij}$ and $\lambda_{ij}$ act identically for each $i$ and $j$. This in turn implies that $[\phi(\rho_{21}),\phi(\rho_{23})]=1$ since $[\rho_{21},\lambda_{23}]=1$. Rewriting the first relation as before we get $[\phi(\rho_{21} \sigma_{13}\rho_{21}), \phi(\sigma_{13})]=1$.

Let $v_2 \in E_2 \smallsetminus \{ 0 \}$. We have $\rho_{21}(v_2) = z_1 + z_2$, with $z_i \in E_i$, and with $z_1 \neq 0$. Thus $z_1$ is a non-zero multiple of $v_1$, since $E_1$ is 1 dimensional.
Arguing as before we now get
\[\langle p_{1,2,3} \big( \rho_{21} \sigma_{13} \rho_{21} (v_2) \big) \rangle  = \langle p_{1,2,3} \big( \rho_{21} \sigma_{13} (z_1) \big) \rangle = U\]

 The group $S_3$ can act on $V_1$ in two ways: via the permutation or the signed permutation representation. In each case however, if $E_1$ is in the $\mu$-eigenspace of $\sigma_{23}$, then $E_2$ is in the $\mu$-eigenspace of $\sigma_{13}$.
Thus ${\sigma_{13}(v_2) = \mu v_2}$ and so $U$ lies in the $\mu$-eigenspace of $\sigma_{13}$, arguing as before. Therefore $U$ also lies in the $\mu$-eigenspace of $\sigma_{12} = \sigma_{23}^{\sigma_{13}}$, and so does $E_3$. This shows that
\[\phi(\rho_{12})\vert_{E_3 \oplus U} = \phi(\rho_{21}^{\sigma_{12}})\vert_{E_3 \oplus U} = \phi(\rho_{21})\vert_{E_3 \oplus U} \]
and therefore that
\[\phi(\rho_{21})\vert_{E_3 \oplus U} = \phi(\lambda_{21})\vert_{E_3 \oplus U} = \phi(\rho_{21} \rho_{12}^{-1} \lambda_{21})\vert_{E_3 \oplus U} = \phi(\epsilon_2 \sigma_{12})\vert_{E_3 \oplus U}\]
But $\epsilon_2 \sigma_{12}(E_3) = E_3$ and so $\rho_{21}$ has a diagram of the form claimed.
\end{proof}

We shall now focus on five- and six-dimensional representations of $\Out{3}$.

\begin{lem}
\label{splitting Delta levels}
\label{lem: gap in levels}
Let $V$ be a $\K$-linear,  six-dimensional \rep of $\Out{3}$, where $\mathrm{char}(\K) \neq 2$. Suppose that, with notation of Definition~\ref{repofW_n}, \[\dim V_1 \oplus V_2 \leqslant 3 .\] Then, if $V_2=\{0\}$, we have a (not necessarily minimal) diagram for $\rho_{21}$ of the form
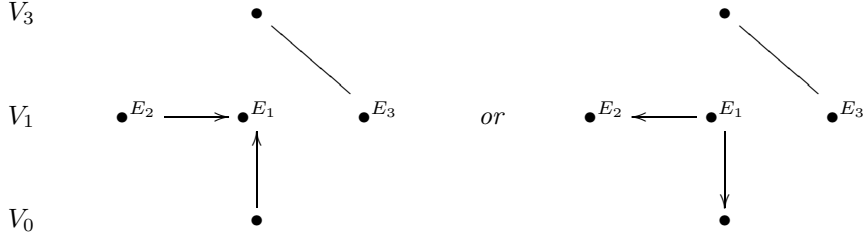
\[ \xymatrix{ V_3 & & \bullet \ar@{-}[dr] & &  & & \bullet \ar@{-}[dr] & & \\
V_1 & \bullet^{E_2} \ar[r] & \bullet^{E_1} \ar@{-}[d] & \bullet^{E_3} & \textrm{or} & \bullet^{E_2} & \bullet^{E_1} \ar[d] \ar[l] & \bullet^{E_3}  \\
V_0 & & \bullet \ar[u] & & & & \bullet  & & }  \]
and if $V_1 = \{0\}$ of the form
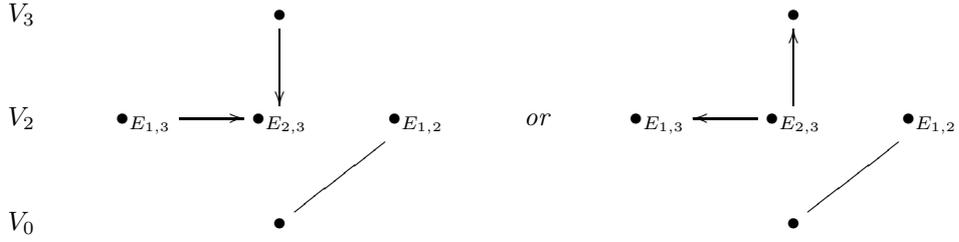
\[ \xymatrix{ V_3 & & \bullet \ar[d] & &  & & \bullet  & & \\
V_2 & \bullet_{E_{1,3}} \ar[r] & \bullet_{E_{2,3}}  & \bullet_{E_{1,2}} & \textrm{or} & \bullet_{E_{1,3}} & \bullet_{E_{2,3}} \ar[u] \ar[l] & \bullet_{E_{1,2}}   \\
V_0 & & \bullet \ar@{-}[ur] & & & & \bullet \ar@{-}[ur]  & & }  \]
In both cases at lest one of $V_0 \oplus V_2$ and $V_1 \oplus V_3$ is $\Out{3}$-invariant.
\begin{proof}
Lemma~\ref{dimensionofV_i} tells us that the dimensions of $V_1$ and $V_2$ are divisible by 3. Hence, by assumption, at least one of $V_1$ and $V_2$ is trivial. If both of them are trivial then Lemma~\ref{diamonds} immediately tells us that the decomposition $V=V_0 \oplus V_3$ is preserved by each $\rho_{ij}$ and $\lambda_{ij}$. Thus the minimal diagram for $\rho_{21}$ is a subdiagram of all the above.

Suppose one of $V_1$, $V_2$ is non-trivial. Without loss of generality let us assume $\dim V_1 \neq 0$. Again by assumption we see that $\dim V_1 = 3$, and hence $\dim E_i = 1$ for all $i$.

Our strategy here is to start with the most general possible diagram for $\rho_{21}$, and then gradually add restriction until we arrive at one of the diagrams described above.

Lemma~\ref{diamonds} allows us to conclude that we have the following diagrams for $\rho_{21}$ and $\rho_{31}$ respectively:
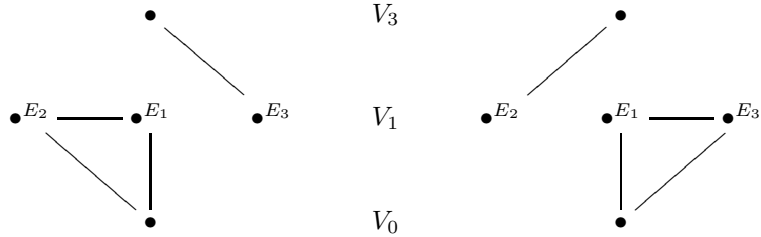
\[ \xymatrix{ & \bullet \ar@{-}[dr] & & V_3 & & \bullet & & \\
\bullet^{E_2} \ar@{-}[r] \ar@{-}[dr] & \bullet^{E_1} \ar@{-}[d] & \bullet^{E_3} & V_1 & \bullet^{E_2} \ar@{-}[ur] & \bullet^{E_1} \ar@{-}[d] \ar@{-}[r] & \bullet^{E_3}  \\
& \bullet & & V_0 & & \bullet \ar@{-}[ur] & & } \]

The element $\Delta$ lies in the centre of $G_3$, and so in particular $[\Delta, \epsilon_1 \sigma_{14} ] =1$. This implies that $\epsilon_1 \sigma_{14}$ preserves the eigenspaces of $\Delta$, which happen to be the direct sums of all subspaces $V_i$ with the index $i$ of a given parity (even for the $(+1)$- and odd for the $(-1)$-eigenspace).
Hence the following is a diagram for $\epsilon_1 \sigma_{14}$:
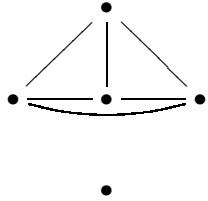
\[ \xymatrix{ & \bullet \ar@{-}[dr] \ar@{-}[d] \ar@{-}[dl] \\
\bullet \ar@{-}[r] \ar@/_/@{-}[rr] & \bullet \ar@{-}[r] & \bullet  \\
& \bullet  } \]
But, in $\Out{3}$, we have $\epsilon_1 \sigma_{14} = \rho_{31} \rho_{21}$, and Example~\ref{ex: component} tells us that
\[ \rho_{31} \Big( p_2 \big( \rho_{21}(V_0) \big) \Big) \leqslant E_2 \oplus V_3 \]
We can therefore conclude that $p_2 \big( \rho_{21}(V_0) \big) = \{0\}$, and so that
we have a diagram for $\rho_{21}$ as follows:
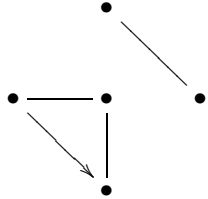
\[ \xymatrix{ & \bullet \ar@{-}[dr]  \\
\bullet \ar@{-}[r] \ar[dr]  & \bullet \ar@{-}[d] & \bullet \\
& \bullet  } \]

Again, by Example~\ref{ex: component}, $\rho_{31}\vert_{E_2 \oplus V_3}$ is an isomorphism. Hence there exists $v \in E_2 \oplus V_3$ such that $\langle \rho_{31} (v) \rangle = E_2$. Since $v \in E_2 \oplus V_3 \leqslant V_1 \oplus V_3$, also $\epsilon_1 \sigma_{14} (v) \in V_1 \oplus V_3$. Now
\[ \epsilon_1 \sigma_{14} (v) = \rho_{21} \rho_{31}(v)\]
and so we conclude that $\rho_{21}$ has a diagram
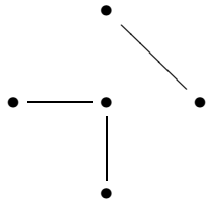
\[ \xymatrix{ & \bullet \ar@{-}[dr]  \\
\bullet \ar@{-}[r] & \bullet \ar@{-}[d] & \bullet \\
& \bullet  } \]

Note that $\rho_{21}$ either has a diagram
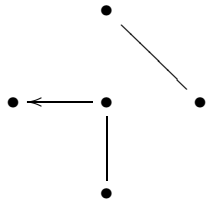
\[ \xymatrix{ & \bullet \ar@{-}[dr]  \\
\bullet  & \bullet \ar[l] \ar@{-}[d] & \bullet \\
& \bullet  } \]
or $p_1 \big( \rho_{21}(E_2) \big) =E_1$, since $\dim E_1 = 1$.

If $\rho_{21}(E_2)$ projects surjectively onto $E_1$, applying $\rho_{31}\rho_{21} = \epsilon_1 \sigma_{14}$ to $E_2$ yields a diagram for $\rho_{31}$ of the form
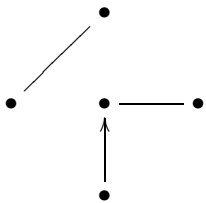
\[ \xymatrix{ & \bullet \ar@{-}[dl]  \\
\bullet  & \bullet \ar@{-}[r]  & \bullet \\
& \bullet \ar[u] } \]
and, after conjugating by $\sigma_{23}$ (see Example~\ref{ex: W_n action}), $\rho_{21}$ has a diagram
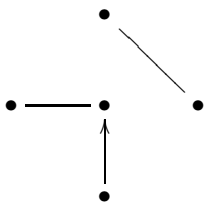
\[ \xymatrix{ & \bullet \ar@{-}[dr]  \\
\bullet  & \bullet \ar@{-}[l]  & \bullet \\
& \bullet \ar[u] } \]
Requiring $\rho_{21}^{\epsilon_1} = \rho_{21}^{-1}$ yields two possibilities for a diagram for $\rho_{21}$:
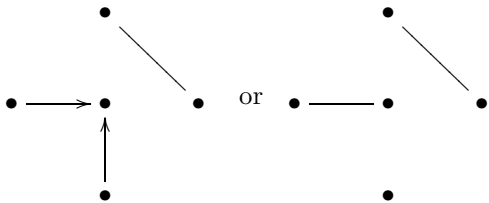
\[ \xymatrix{ & \bullet \ar@{-}[dr]  \\
\bullet \ar[r]  & \bullet   & \bullet \\
& \bullet \ar[u] } \ \xymatrix{ \\ \textrm{ or }}
\xymatrix{ & \bullet \ar@{-}[dr]  \\
\bullet  & \bullet \ar@{-}[l]  & \bullet \\
& \bullet  } \]
The first diagram is as required. The second diagram gives a required diagram after applying Lemma~\ref{lem: annoying case}.

We still have to consider the case of a diagram
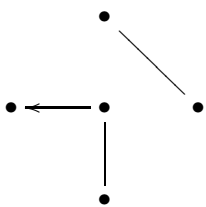
\[ \xymatrix{ & \bullet \ar@{-}[dr]  \\
\bullet  & \bullet \ar[l] \ar@{-}[d] & \bullet \\
& \bullet  } \]
for $\rho_{21}$. Applying $\epsilon_1 \sigma_{14} = \rho_{21} \rho_{31}$ to $V_0$ yields a diagram for $\rho_{21}$ of the form
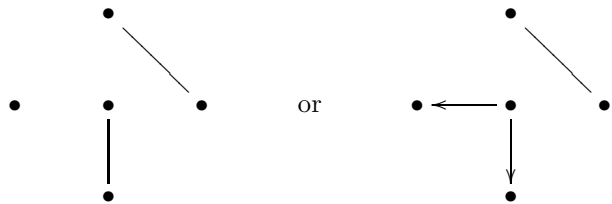
\[ \xymatrix{ & \bullet \ar@{-}[dr] & &  & & \bullet \ar@{-}[dr] & & \\
\bullet  & \bullet \ar@{-}[d] & \bullet & \textrm{ or } & \bullet & \bullet \ar[d] \ar[l] & \bullet  \\
& \bullet & & & & \bullet  & & } \]
The second of these diagrams is as required.

Let us now focus on the first of the above diagrams. Note that, by Example~\ref{ex: W_n action}, this is also a diagram for $\lambda_{21}$, and that a diagram for $\rho_{12}^{-1}$ is as follows:
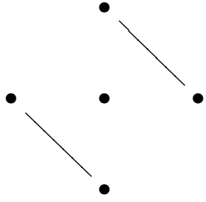
\[ \xymatrix{ & \bullet \ar@{-}[dr]  \\
\bullet \ar@{-}[dr] & \bullet  & \bullet   \\
& \bullet  } \]
Let $v_1$ be a generator of $E_1$. Apply $\epsilon_2 \sigma_{12} = \rho_{21} \rho_{12}^{-1} \lambda_{21}$ to $v_1$ and observe that $\epsilon_2 \sigma_{12} (v_1) = v_2$, a generator of $E_2$. Now let $x$ be the $E_1$ component of $\lambda_{21}(v_1)$. Note that $\rho_{12}^{-1} \lambda_{21}(v_1)$ has a non-trivial $E_1$ component \iff $x$ is not zero. But such a non-trivial component yields a non-zero component in $E_1 \oplus V_0$ of $\rho_{21} \rho_{12}^{-1} \lambda_{21}(v_1)$. This is impossible, since $\epsilon_2 \sigma_{12} (v_1) = v_2$ has no such components. Thus $x=0$, $\lambda_{21}(v_1)$ lies in $V_0$, and
\[\rho_{12}^{-1}\vert_U :U  \to E_2\]
is an isomorphism, where $U=\langle \lambda_{21} (v_1) \rangle$. Hence $\rho_{12}\vert_{E_2} : E_2 \to U$ is an isomorphism as well.

We claim that $\rho_{ij}^{\pm 1}\vert_{E_j}, \lambda_{ij}^{\pm 1}\vert_{E_j} : E_j \to U$ are all isomorphisms. We have established this for $\lambda_{21}$ and $\rho_{12}$. Conjugating by $\epsilon_1$ and $\epsilon_2$ establishes the claim also for $\rho_{21}^{-1}, \rho_{21}, \lambda_{21}^{-1}, \rho_{12}^{-1}, \lambda_{12}$ and $\lambda_{12}^{-1}$. Using the fact that $\epsilon_1 \sigma_{14} = \rho_{31} \rho_{21}$ preserves $V_1 \oplus V_3$ we immediately conclude that the claim also holds for $\rho_{31}^{-1}$, and hence in particular also for $\rho_{13}$ (repeating the argument above). Now the relation $\epsilon_3 \sigma_{34} = \rho_{13} \rho_{23}$ establishes the claim for $\rho_{23}$, and the claim follows.

Our calculations enable us to deduce that diagrams for $\rho_{21}$ and $\lambda_{23}$ respectively are as follows
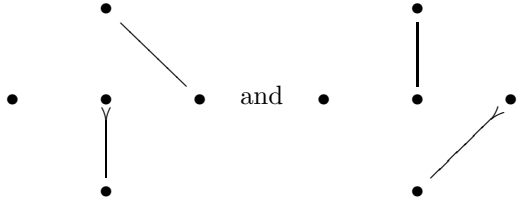
\[ \xymatrix{ & \bullet \ar@{-}[dr]  \\
\bullet  & \bullet \ar@{ >-}[d] & \bullet   \\
& \bullet  } \ \xymatrix{ \\ \textrm{ and }} \
\xymatrix{ & \bullet \ar@{-}[d]  \\
\bullet  & \bullet  & \bullet  \ar@{ >-}[dl] \\
& \bullet }
\]
But $\rho_{21}$ and $\lambda_{23}$ commute, and this together with the fact that $\rho_{21}(E_1) = U = \lambda_{23}^{-1}(E_3)$ yields a diagram for $\rho_{21}$ of the form
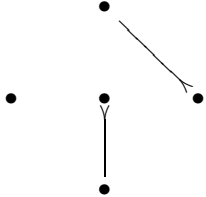
\[ \xymatrix{ & \bullet   \\
\bullet  & \bullet \ar@{ >-}[d] & \bullet \ar@{ >-}[ul] \\
& \bullet  } \]
In particular Example~\ref{ex: away arrow} implies that $\dim V_0 \neq 0$.

Now let us define $A_2 = \rho_{21}(E_3) \leqslant V_3$. Note that $\dim A_2 = 1$. Since $\rho_{21}$ commutes with $\lambda_{23}$, examining the respective diagrams yields $\lambda_{23} (A_2) = E_1$. Now, observing that each $\epsilon_i$ preserves each subspace of $E_I$, we see that in fact for all $i$
\[A_2 = \rho_{2i}^{\pm1}(E_j)  = \lambda_{2i}^{\pm1}(E_j) \]
where $j$ satisfies $\{i,j\} = \{1,3\}$. We can define $A_1$ and $A_3$ similarly.

The relations $[\rho_{21},\rho_{31}]=[\rho_{23},\rho_{13}]=1$, together with the structure of our diagrams, tell us that $A_2 \cap ( A_1 +  A_3 )=\{0\}$. The relation $[\rho_{32},\rho_{12}]=1$ informs us that $A_1 \cap A_3 = \{ 0\}$ and so that $\dim ( A_1 \oplus A_2 \oplus A_3) =3$. This is a contradiction, since $V_0 \neq \{0\}$ and so $\dim V_3 = 2$.

We have thus shown that $\rho_{21}$ has a diagram as claimed. Observe that, since the subgroup $W_3 < \Out 3$ preserves each $V_i$ by construction, having a diagram for $\rho_{21}$ of the form described in the statement of this lemma immediately implies that at least one of $V_0 \oplus V_2$ and $V_1 \oplus V_3$ is preserved by $\Out{3} = \langle W_3, \rho_{21} \rangle$.
\end{proof}
\end{lem}

\begin{lem}
\label{dimension 6 bad case}
\label{lem: dim 6 new case}
Let $V$ be a $\K$-linear,  six-dimensional \rep of $\Out{3}$, where $\K$ is a field of characteristic other than 2 or 3. Suppose that, with notation of Definition~\ref{repofW_n}, $\dim V_1 \oplus V_2 = 6$. Then $V = V_1 \oplus V_2$ as an $\Out{3}$-module.
\begin{proof}
If $\dim V_1 = 6$ or $\dim V_2 = 6$ then the result is trivial.

Suppose that $\dim V_1 = \dim V_2 = 3$ and so $V = V_1 \oplus V_2$ as a vector space.
We know (using Maschke's Theorem and our assumption on $\mathrm{char}(\K)$) that each $V_i$ (for $i = 1,2$) is either a sum of standard and trivial or a sum of signed standard and determinant \reps of $S_3$; we can therefore pick vectors $v_i \in E_i, w_i \in E_{\{1,2,3\} \smallsetminus \{i\}}$ so that each $v_i - v_j$ and $w_i - w_j$ is an eigenvector of an element of $S_3 \smallsetminus \{1\}$.

We have a diagram for $\rho_{21}$ of the form
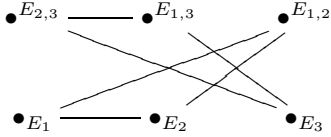
\[ \xymatrix{ \bullet^{E_{2,3}} \ar@{-}[r]  \ar@{-}[drr] & \bullet^{E_{1,3}} \ar@{-}[dr] & \bullet^{E_{1,2}}  \\
\bullet_{E_{1}} \ar@{-}[r] \ar@{-}[urr] & \bullet_{E_{2}} \ar@{-}[ur] & \bullet_{E_{3}} } \]
and analogously one for $\rho_{31}$ of the form
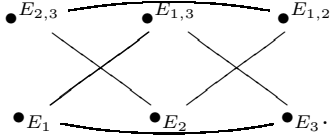
\[ \xymatrix{ \bullet^{E_{2,3}} \ar@{-}[dr]  \ar@/^/@{-}[rr] & \bullet^{E_{1,3}} \ar@{-}[dl] & \bullet^{E_{1,2}} \ar@{-}[dl]  \\
\bullet_{E_{1}} \ar@{-}[ur] \ar@/_/@{-}[rr] & \bullet_{E_{2}} & \bullet_{E_{3}} . \ar@{-}[ul] } \]
Since $S_4$ commutes with $\Delta$, its action has to preserve the $(+1)$-eigenspace of $\Delta$ (which is equal to $V_1$ in our case) as well as the $(-1)$-eigenspace (which equals $V_2$ in this case).
We also have $[\epsilon_1 , \Delta] =1$, and so
$\epsilon_1 \sigma_{14}  = \rho_{31} \rho_{21}$ preserves $V_2$. Hence, evaluating $\rho_{31} \rho_{21}$ on $E_{1,2}$ (an observing that $\dim E_I \leqslant 1$ for all $I$) gives us either a diagram for $\rho_{21}$ of the form
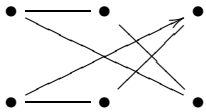
\[ \xymatrix{ \bullet \ar@{-}[r]  \ar@{-}[drr] & \bullet \ar@{-}[dr] & \bullet  \\
\bullet \ar@{-}[r] \ar[urr] & \bullet \ar@{-}[ur] & \bullet } \]
or a diagram for $\rho_{31}$ of the form
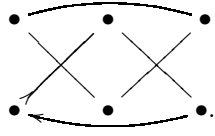
\[ \xymatrix{ \bullet \ar@{-}[dr]  \ar@/^/@{-}[rr] & \bullet \ar@{-}[dl] & \bullet \ar@{-}[dl]  \\
\bullet \ar@{ >-}[ur] & \bullet & \bullet . \ar@/^/[ll] \ar@{-}[ul]  } \]

Suppose (for a contradiction) that we are in the latter case. Evaluating $\rho_{21} \rho_{31}$ on $E_1$ (and observing that the diagrams for $\rho_{31}$ and $\rho_{21}$ are related by conjugation by $\sigma_{23}$) yields diagrams for $\rho_{21}$ and $\rho_{31}$ respectively of the form
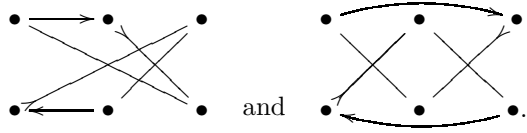
\[ \xymatrix{ \bullet \ar[r]  \ar@{-}[drr] & \bullet \ar@{ >-}[dr] & \bullet  \\
\bullet  \ar@{ >-}[urr] & \bullet \ar[l] \ar@{-}[ur] & \bullet } \ \xymatrix{ \\ \textrm{ and } \\ } \
\xymatrix{ \bullet \ar@{-}[dr]  \ar@/^/[rr] & \bullet \ar@{-}[dl] & \bullet \ar@{ >-}[dl]  \\
\bullet \ar@{ >-}[ur] & \bullet & \bullet . \ar@/^/[ll] \ar@{-}[ul]  }
 \]
Now $\rho_{21} \rho_{31} (E_1) = E_3$ and $\rho_{31} \rho_{21} (E_1) = E_2$. But $\rho_{31}$ commutes with $\rho_{21}$, which yields a contradiction.

We can repeat the argument after evaluating $\rho_{31} \rho_{21}$ on $E_3$ and conclude that we have a diagram for $\rho_{21}$ of the form
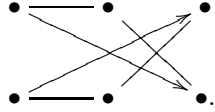
\[ \xymatrix{ \bullet \ar@{-}[r]  \ar[drr] & \bullet \ar@{-}[dr] & \bullet  \\
\bullet \ar@{-}[r] \ar[urr] & \bullet \ar@{-}[ur] & \bullet .} \]

Two diagram chases, starting at $E_3$ and $E_{1,2}$, show $\rho_{21}^{-1} = \rho_{21}^{\epsilon_1}$ requires $\rho_{21}$ to have a diagram of the form
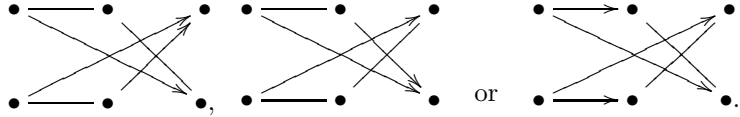
\[ \xymatrix{ \bullet \ar@{-}[r]  \ar[drr] & \bullet \ar@{-}[dr] & \bullet  \\
\bullet \ar@{-}[r] \ar[urr] & \bullet \ar[ur] & \bullet , } \
\xymatrix{ \bullet \ar@{-}[r]  \ar[drr] & \bullet \ar[dr] & \bullet  \\
\bullet \ar@{-}[r] \ar[urr] & \bullet \ar@{-}[ur] & \bullet }
\ \xymatrix{ \\ \textrm{ or }} \
\xymatrix{ \bullet \ar[r]  \ar[drr] & \bullet \ar@{-}[dr] & \bullet  \\
\bullet \ar[r] \ar[urr] & \bullet \ar@{-}[ur] & \bullet .}
\]

Suppose we are in the third case. We have diagrams for $\rho_{21}$ and $\rho_{31}$ respectively
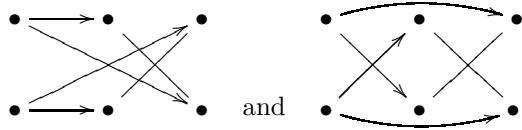
\[ \xymatrix{ \bullet \ar[r]  \ar[drr] & \bullet \ar@{-}[dr] & \bullet  \\
\bullet \ar[r] \ar[urr] & \bullet \ar@{-}[ur] & \bullet } \ \xymatrix{ \\ \textrm{ and }} \
\xymatrix{ \bullet \ar[dr]  \ar@/^/[rr] & \bullet \ar@{-}[dl] & \bullet \ar@{-}[dl]  \\
\bullet \ar[ur] \ar@/_/[rr] & \bullet & \bullet .  \ar@{-}[ul]  }
\]
Evaluating $\epsilon_1 \sigma_{14} = \rho_{31} \rho_{21}$ on $E_1$ (and observing that $\epsilon_1 \sigma_{14}(V_2) = V_2$)  yields a diagram for $\rho_{21}$ of the form
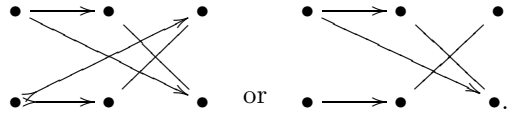
\[ \xymatrix{ \bullet \ar[r]  \ar[drr] & \bullet \ar@{-}[dr] & \bullet  \\
\bullet \ar@{ >->}[r] \ar@{ >->}[urr] & \bullet \ar@{-}[ur] & \bullet } \ \xymatrix{ \\ \textrm{ or }} \
\xymatrix{ \bullet \ar[r]  \ar[drr] & \bullet \ar@{-}[dr] & \bullet  \\
\bullet \ar[r]  & \bullet \ar@{-}[ur] & \bullet . } \]
The first case is impossible, since we could have
\[E_1 = 1(E_1) = \epsilon_1 \rho_{21} \epsilon_1 \rho_{21} (E_1) \leqslant E_2 \oplus E_{1,2} .\]
After repeating the argument for $E_{2,3}$ we conclude that we have diagrams for $\rho_{21}$ and $\rho_{31}$ respectively as follows
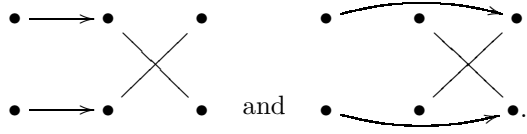
\[ \xymatrix{ \bullet \ar[r]   & \bullet  & \bullet \ar@{-}[ld]  \\
\bullet \ar[r]  & \bullet  & \bullet  \ar@{-}[ul]} \ \xymatrix{ \\ \textrm{ and }} \
\xymatrix{ \bullet   \ar@/^/[rr] & \bullet \ar@{-}[dr] & \bullet   \\
\bullet  \ar@/_/[rr] & \bullet \ar@{-}[ur] & \bullet  .   }
\]

Suppose that $\sigma_{14}$ preserves each $E_i$. Then so does $\sigma_{24} = \sigma_{14}^{\sigma_{12}}$. But ${\sigma_{24} = \sigma_{12}^{\sigma_{14}}}$, and $\sigma_{12}(E_1) = E_2$. This is a contradiction. We can apply an analogous argument to the $\sigma_{14}$-action on the subspaces $E_{i,j}$. Now we easily deduce from $\epsilon_1 \sigma_{14} = \rho_{31} \rho_{21}$ that \[\rho_{21}(E_1) \not\leqslant E_1 \textrm{ and } \rho_{21}(E_{2,3}) \not\leqslant E_{2,3} .\]

We can now evaluate $\epsilon_1 \sigma_{14} = \rho_{31} \rho_{21}$ and $\rho_{21}^{-1} = \rho_{21}^{\epsilon_1}$ on $E_1$ and $E_{2,3}$ and conclude that we have a diagram for $\rho_{21}$ of the form
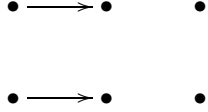
\[ \xymatrix{ \bullet \ar[r]   & \bullet  & \bullet  \\
\bullet \ar[r]  & \bullet  & \bullet}
\]
which shows that both $V_1$ and $V_2$ are $\Out{n}$-invariant.

Suppose now that we are in one of the first two cases, namely that there is a diagram for $\rho_{21}$ of the form
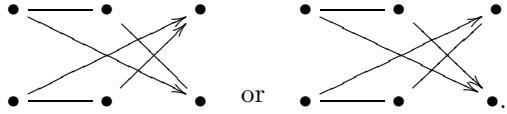
\[ \xymatrix{ \bullet \ar@{-}[r]  \ar[drr] & \bullet \ar@{-}[dr] & \bullet  \\
\bullet \ar@{-}[r] \ar[urr] & \bullet \ar[ur] & \bullet } \ \xymatrix{ \\ \textrm{ or }} \
\xymatrix{ \bullet \ar@{-}[r]  \ar[drr] & \bullet \ar[dr] & \bullet  \\
\bullet \ar@{-}[r] \ar[urr] & \bullet \ar@{-}[ur] & \bullet. } \]
Verifying that $\rho_{31} \rho_{21}$ keeps $V_1$ and $V_2$ invariant immediately tells us that in fact we have a diagram for $\rho_{21}$ of the form
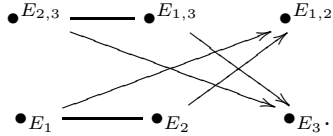
\[ \xymatrix{ \bullet^{E_{2,3}} \ar@{-}[r]  \ar[drr] & \bullet^{E_{1,3}} \ar[dr] & \bullet^{E_{1,2}}  \\
\bullet_{E_1} \ar@{-}[r] \ar[urr] & \bullet_{E_2} \ar[ur] & \bullet_{E_3} .} \]
The element $\rho_{31}$ keeps $E_2$ and $E_{1,3}$ invariant, and so, observing that $\epsilon_1 \sigma_{14} = \rho_{21} \rho_{31}$ preserves $V_1 \oplus V_3$, we actually have diagrams for $\rho_{21}$ and $\rho_{31}$ respectively
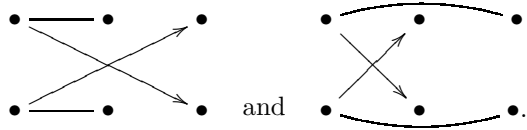
\[ \xymatrix{ \bullet \ar@{-}[r]  \ar[drr] & \bullet  & \bullet  \\
\bullet \ar@{-}[r] \ar[urr] & \bullet  & \bullet } \ \xymatrix{ \\ \textrm{ and }} \
\xymatrix{ \bullet   \ar@/^/@{-}[rr] \ar[dr] & \bullet  & \bullet   \\
\bullet  \ar@/_/@{-}[rr] \ar[ur] & \bullet  & \bullet  .   }
\]
But, in order for $\rho_{21} \rho_{31}$ to keep $V_1$ and $V_2$ invariant, we need to have
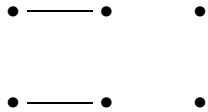
\[ \xymatrix{ \bullet \ar@{-}[r]   & \bullet  & \bullet  \\
\bullet \ar@{-}[r]  & \bullet  & \bullet } \]
as a diagram for $\rho_{21}$. This finishes the proof.
\end{proof}
\end{lem}

Now let us investigate 5-dimensional representations of $\Out{3}$ -- we hope to be able to say more in this case!

\begin{prop}
\label{reps prop}
Let $V$ be a 5-dimensional, $\K$-linear \rep of $\Out{3}$, where $\K$ is a field of characteristic other than 2 or 3. Suppose that, with the notation of Definition~\ref{repofW_n}, $V \neq V_0 \oplus V_3$. Then $V=V_0 \oplus V_1 \oplus V_2 \oplus V_3$ is a decomposition of $\Out{3}$-modules, and, as $S_4$-modules, $V_0$ is a sum of trivial, $V_1$ of standard, $V_2$ of signed standard, and $V_3$ of determinant representations.
\begin{proof}
Since $\dim V =5$, we have $V_1 = \{ 0 \}$ or $V_2 = \{ 0 \}$. Let us suppose that we have the latter, the other case being entirely similar.

\medskip
\noindent \textbf{Step 0:} We first claim that $V_0$ is a sum of trivial $S_4$-modules.
\medskip

Lemma~\ref{splitting Delta levels} gives us two possibilities for a diagram for $\rho_{21}$, namely
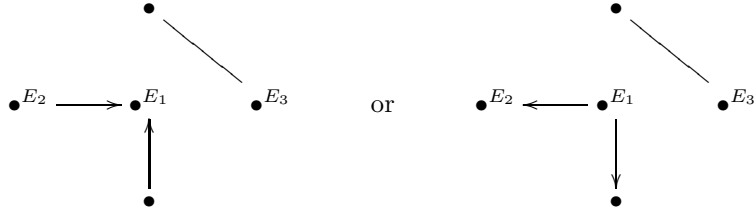
\[ \xymatrix{ & \bullet \ar@{-}[dr] & &  & & \bullet \ar@{-}[dr] & & \\
 \bullet^{E_2} \ar[r] & \bullet^{E_1} \ar@{-}[d] & \bullet^{E_3} & \textrm{or} & \bullet^{E_2} & \bullet^{E_1} \ar[d] \ar[l] & \bullet^{E_3}   \\
 & \bullet \ar[u] & & & & \bullet  & & }  \]
The same lemma also tells us that $V/(V_1 \oplus V_3)$ is a \rep of $\Out{3}$. Its dimension is at most 2 and therefore Lemma~\ref{killV_4} tells us that it is a direct sum of two trivial \reps of $\Out{3}$ (since we know how $\epsilon_1$ acts), and so the same statement holds for $V / (V_1 \oplus V_3)$ as an $S_4$-module. Hence it also holds for $V_0$, since $V_0$ is an $S_4$-module isomorphic to $V/(V_1 \oplus V_3)$.

Note that an identical argument shows that $V_3$ is a sum of determinant $S_4$-modules in the case when $V_1 = \{0\}$.

\medskip
\noindent \textbf{Step 1:} We now claim that $V_0 \oplus V_1$ is $\Out{3}$-invariant.
Suppose for a contradiction that it is not the case.
\medskip

Let $U$ be the projection of  $\rho_{21}(E_3)$ onto $V_3$. Note that $\dim U = 1$ since we have assumed $V_0\oplus V_1$ not to be $\Out{3}$-invariant. Our aim now is to show that $U$ is $\Out{3}$-invariant.

If $V_3$ is $\Out{3}$-invariant, then it is an $\Out{3}$-module of dimension at most two, and hence we can use Lemma~\ref{repsofdim2} to conclude that it is in fact a sum of determinant representations. Hence, in particular, $U$ is $\Out 3$-invariant.

Now suppose that $V_3$ is not $\Out{3}$-invariant.
Checking that $\rho_{31} \rho_{21}(E_3) = \rho_{21} \rho_{31}(E_3)$ on both of our diagrams for $\rho_{21}$ yields that $\rho_{31}(U) = U$.
Note that $U$ is the unique non-trivial invariant subspace of $V_3$ for both $\rho_{21}$ and $\rho_{31}$, as otherwise $V_3$ would be invariant under the action of \[\langle S_3, \rho_{21} \rangle = \langle S_3, \rho_{31} \rangle = \Out 3 .\] Hence $U$ is $\sigma_{23}$-invariant. But $V_3$ is a 2-dimensional $S_3$-module, and if it were irreducible, then the trace of each transposition would be zero. Hence $V_3$ is a sum of two 1-dimensional modules of $S_3$, and therefore $S_3$ preserves $U$. From this we conclude that $\Out 3$ preserves $U$.

Lemma~\ref{lem: 2-dim reps} informs us that $U$ is a determinant \rep of $\Out 3$. Since $\rho_{21}^{-1} = \rho_{21}^{\epsilon_1}$, we must have
\[ \forall v \in E_3 \, : \, \rho_{21}(v) \in v + U .\] Using similar relations we establish that, when restricted to $E_3 \oplus U$, $\lambda_{21}$ acts as $\rho_{21}$, and $\rho_{12}$ acts as $\rho_{21}^{\pm1}$. Hence, taking $v \in E_3$,
\[  v + (2 \mp1) u = \rho_{21}^{2 \mp 1}(v)= \rho_{21} \rho_{12}^{-1} \lambda_{21}(v)=\epsilon_1 \sigma_{12} ( v) \in E_3\]
where $u = \rho_{21}(v) -v \in U$. This shows that $u=0$, and hence $V_0\oplus V_1$ is $\Out{3}$-invariant, which is the desired contradiction.

We have thus shown that there is a diagram for $\rho_{21}$ of the form
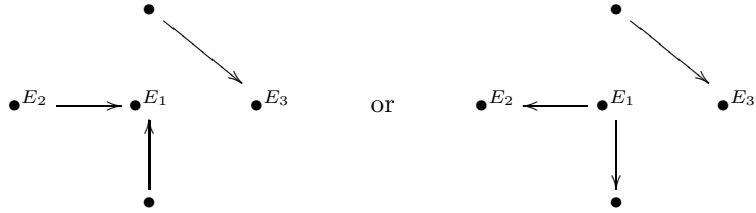
\[ \xymatrix{ & \bullet \ar[dr] & &  & & \bullet \ar[dr] & & \\
 \bullet^{E_2} \ar[r] & \bullet^{E_1} \ar@{-}[d] & \bullet^{E_3} & \textrm{or} & \bullet^{E_2} & \bullet^{E_1} \ar[d] \ar[l] & \bullet^{E_3}   \\
 & \bullet \ar[u] & & & & \bullet  & & }  \]

\medskip
\noindent \textbf{Step 2:} We claim that $V_1$ is a standard $S_4$-module.
\medskip

As an $S_3$-module, both $V_1$ and $V_2$ are sums of one standard and one either trivial or determinant representation. The branching rule tells us that therefore, as $S_4$-modules, each of the subspaces can be either a standard or a signed standard representation, or the one corresponding to partition $(2,2)$. The last case is ruled out by Lemma~\ref{killV_4}, since $(V_0 \oplus V_1)/V_0$ is clearly not a sum of trivial and determinant $\Out{3}$-modules.

Focusing only on $V_1$, we have a diagram for $\rho_{21}$ of the form
\[ \xymatrix{
 \bullet \ar[r] & \bullet  & \bullet & \textrm{or} & \bullet & \bullet \ar[l] & \bullet   }  \]
Note that in both cases these are the minimal diagrams for $\rho_{21}$ when restricted to $V_1$, since otherwise $\sigma_{12}$ could not permute $E_1$ and $E_2$.

Let us pick vectors $v_i \in E_i$ in such a way that each $v_i-v_j$ is an eigenvector of $\sigma_{ij}$. Let us also set $v = v_1 + v_2 + v_3$. The way in which $S_4$ acts on $V_1$ in our case is determined by one parameter; we can calculate it by finding $\mu \in \C$ such that $v_1 + \mu v$ is an eigenvector of $\sigma_{14}$. The eigenvalue of this eigenvector will also determine the way in which $S_4$ acts. Let us note that we can also find this parameter $\mu$ by computing $\sigma_{14}(v_2 - v_1) = \mu v + v_2$.

In the case of the first diagram for $\rho_{21}$, we immediately see that
\[\sigma_{14}(v_2-v_1) = \epsilon_{1} \rho_{21} \rho_{31} (v_2-v_1) \in E_1 \oplus E_2 ,\]
and hence $\mu = 0$. Now both $\rho_{23}$ and $\rho_{31}$ preserve $E_1$, and so observing that $\rho_{21}^{-1} = [\rho_{23}^{-1} ,\rho_{31}^{-1}]$ yields that $\rho_{21}$ acts trivially on $E_1$.
By an analogous argument so does $\rho_{31}$. Hence $\sigma_{14}(v_1) = \epsilon_1(v_1)$. In our case this shows that we are dealing with a standard representation; if however $V_1$ is trivial, $\epsilon_1$ acts as plus one on the appropriate vector, and wee see that $V_2$ is a signed standard $S_4$-module.

In the case of the second diagram we immediately see two eigenspaces of $\sigma_{14}$, namely $E_2$ and $E_3$. These spaces are interchanged by the action of $\sigma_{23}$ which commutes with $\sigma_{14}$, and hence must have the same eigenvalue. In a standard or a signed standard \rep of $S_4$ each $\sigma_{ij}$ has always exactly two repeated eigenvalues, and it is this eigenvalue that determines the representation. It is enough for us then to find a third eigenvector of $\sigma_{14}$ and compute its eigenvalue. The vector must have a non-trivial $E_1$-component, and our diagram tells us that it is enough to check how $\sigma_{14}$ acts on $E_1$. By an argument similar to the one above we show that $\epsilon_1 \sigma_{14}(v_1) \in v_1 + E_2 \oplus E_3$, and the claim follows.

\medskip
\noindent \textbf{Step 3:} We now claim that $V_3$ is a sum of determinant $S_4$-representations.
\medskip

As an $S_3$-module, $V_1 \oplus V_3$ is a sum of one standard, one trivial and some number of determinant representations (depending on the dimension of $V_3$). We have already found one standard \rep of $S_4$, and the branching rule tells us that there can only be determinant \reps of $S_4$ left. If at least one of them does not lie entirely in $V_3$, then it would appear in $(V_1 \oplus V_3) / V_3$ by Schur's Lemma. This is not possible, since $(V_1 \oplus V_3) / V_3$ is a standard \rep of $S_4$. Hence all the other irreducible $S_4$-modules lie within $V_3$.

\medskip
\noindent \textbf{Step 4:} Our last claim is that each $V_i$ is $\Out{3}$-invariant.
\medskip

We have already shown this for $V_1\oplus V_0$. We have just shown that $V_3$ is $S_4$ invariant, and so, $\rho_{21} \rho_{31} = \epsilon_1 \sigma_{14}$ keeping $V_3$ invariant yields a diagram for $\rho_{21}$ of the form
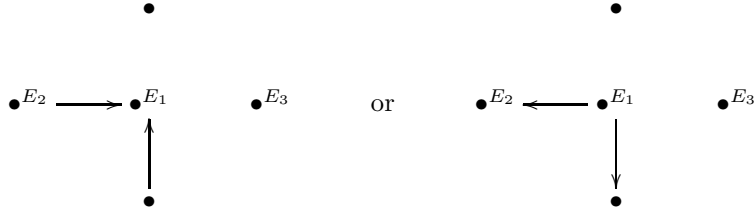
\[ \xymatrix{ & \bullet  & &  & & \bullet  & & \\
 \bullet^{E_2} \ar[r] & \bullet^{E_1} \ar@{-}[d] & \bullet^{E_3} & \textrm{or} & \bullet^{E_2} & \bullet^{E_1} \ar[d] \ar[l] & \bullet^{E_3}   \\
 & \bullet \ar[u] & & & & \bullet  & & }  \]
We have already shown in step 2 that in both cases $\rho_{21}(v) \in v + E_2 \oplus V_0$ for each $v \in E_1$. Also, $(V_0\oplus V_1) / V_1$ is an $\Out{3}$-module of dimension at most 2, and hence is described by Lemma~\ref{repsofdim2}. In particular $\rho_{12}(w) = w + E_1$ for all $w \in V_0$. Analogous statements hold for $\rho_{31}$ and so observing that $\sigma_{23}$ acts as $\pm 1$ on $E_1 \oplus V_0$ and that $\rho_{21} = \rho_{31}^{\sigma_{23}}$ yields that $\epsilon_1 \sigma_{12} = \rho_{31} \rho_{21} (V_0)$ has a non-trivial $V_1$-component \iff $\rho_{21}(V_0)$ does, and similarly that $\epsilon_1 \sigma_{12} = \rho_{31} \rho_{21} (V_1)$ has a non-trivial $V_0$-component \iff $\rho_{21}(V_0)$ does. Hence
 we have a diagram
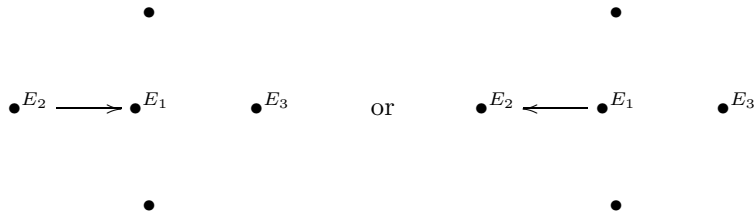
\[ \xymatrix{ & \bullet  & &  & & \bullet  & & \\
 \bullet^{E_2} \ar[r] & \bullet^{E_1} & \bullet^{E_3} & \textrm{or} & \bullet^{E_2} & \bullet^{E_1}  \ar[l] & \bullet^{E_3}   \\
 & \bullet  & & & & \bullet  & & }  \]
 for $\rho_{21}$, which was what we claimed.
\end{proof}
\end{prop}

\begin{lem}
\label{lem: commuting delta}
Suppose $V$ is a $\K$-linear \rep of $\Out{3}$, such that, using notation of Definition~\ref{repofW_n}, $V_0 \oplus V_2$ and $V_1 \oplus V_3$ are $\Out 3$-invariant. Then the representation factors through the natural projection \[\pi: \Out{3} \to \GL_3(\Z) .\]
\begin{proof}
Note that $\phi(\Delta)$ lies in the product \[Z(\GL(V_0 \oplus V_2)) \times Z(\GL(V_1 \oplus V_3)) \] of the centres of the general linear groups of the components $V_0 \oplus V_2$ and $V_1 \oplus V_3$. Therefore we have $\phi(\rho_{ij}) = \phi(\rho_{ij})^{\phi(\Delta)} = \phi(\lambda_{ij})$ for each $i \neq j$, and so $\phi$ factors as
\[ \xymatrix{ \Out{3} \ar[r]^\phi \ar[d]& \GL(V) \\
 \Out{3}/ \langle \! \langle \{ \rho_{ij}\lambda_{ij}^{-1} \mid i\neq j \} \rangle \! \rangle \ar[r]^{\quad \quad \quad \quad \cong} & \GL_3(\Z)   \ar[u] }  \]
This finishes the proof.
\end{proof}
\end{lem}

Observe an immediate consequence of the above.
\begin{lem}
\label{lem: 5-dim reps}
Suppose $V$ is a $\K$-linear \rep of $\Out{3}$ of dimension at most 5, where the characteristic of $\K$ is not 2 or 3. Then the representation factors through the natural projection $\pi: \Out{3} \to \GL_3(\Z)$.
\begin{proof}
Using the notation of Definition~\ref{repofW_n}, we have $V= V_0 \oplus V_1 \oplus V_2 \oplus V_3$ as a vector space.
Suppose first that $V_1 \oplus V_2$ is trivial. Then Lemma~\ref{diamonds} tells us that $V = V_0 \oplus V_3$ as an $\Out{3}$-module.

Supposing that $V_1 \oplus V_2 \neq \{0\}$ allows us to use Proposition~\ref{reps prop}, and conclude that each $V_i$ is $\Out 3$-invariant. We can now use Lemma~\ref{lem: commuting delta} and finish the proof.
\end{proof}
\end{lem}

Before proceeding further we need to recall a standard fact of representations theory.

\begin{prop}
\label{prop: irreducibility of K-module}
Let $A$ be the kernel of the map $\SL_3(\Z) \to \SL_3(\Z_2)$ induced by the surjection $\Z \to \Z_2$. Let $V$ be the standard, 3-dimensional $\K$-linear representation of $\GL_3(\Z)$. Suppose further that $\K$ is a field of characteristic 0 or at least 3. Then $U=\mathrm{Sym}^2(V^*)$, the second symmetric power of the dual module of $V,$ is irreducible as an $A$-module.
\begin{proof}
Let $U \leqslant V$ be an irreducible $A$-submodule of $V$, and let $\{ v_1, v_2, v_3 \}$ be the standard basis of $V$. Suppose $v \in U \smallsetminus \{0\}$. Then
\[v = \sum_{i \leqslant j} \mu_{i j} \, v_i \otimes v_j \]
for some collection of scalars $\mu_{ij}$.

We are going to abuse notation by using the symbols $\epsilon_i$ and $\rho_{ij}$ to denote the images of respective elements under $\pi : \Out 3 \to \GL_3(\Z)$. Note that $\epsilon_i \epsilon_j \in A$ and $\rho_{ij}^2 \in A$ for each appropriate $i\neq j$. Now
\[ \epsilon_1 \epsilon_2 (v) -v = -2 \, \mu_{23} \, v_2 \otimes v_3 -2 \, \mu_{13} \, v_1 \otimes v_3 \]
and hence
\[ \epsilon_1 \epsilon_3\big(\epsilon_1 \epsilon_2 (v) -v\big) - v = 4 \, \mu_{13} \, v_1 \otimes v_3.\]
Hence, if $\mu_{ij} \neq 0$ for some $i \neq j$, then $v_i \otimes v_j \in U$.

Furthermore
\[ \rho_{13}^2 (v_1 \otimes v_3) - v_1 \otimes v_3 = -2 \, v_1 \otimes v_1 \]
and
\[ \rho_{23}^2 (v_1 \otimes v_3) - v_1 \otimes v_3 = -2 \, v_1 \otimes v_2 \]
and therefore if $\mu_{ij} \neq 0$ for some $i \neq j$, then $U=V$.

Suppose that
\[ v = \sum_{i } \mu_{i i} \, v_i \otimes v_i \]
\Wlog let us assume that $\mu_{1 1} \neq 0$.
Then
\[ \rho_{21}^2 (v) - v = -\mu_{1 1} \big(2 \, v_1 \otimes v_2 -4 \, v_2 \otimes v_2 \big) = v' \in U \]
We can now apply our argument to $v'$ and conclude that $U=V$.
\end{proof}
\end{prop}

We are now ready for the main result of this section.

\begin{thm}
\label{result A}
\label{thm: rep thm}
Suppose $V$ is a $\K$-linear \rep of $\Out{3}$ of dimension at most 6, where the characteristic of $\K$ is not 2 or 3. Then the representation factors through the natural projection $\pi: \Out{3} \to \GL_3(\Z)$.
\begin{proof}
Let $\phi : \Out 3 \to \GL(V)$ be our representation. Using the notation of Definition~\ref{repofW_n}, we have $V= V_0 \oplus V_1 \oplus V_2 \oplus V_3$ as a vector space.
We need to consider a number of cases.

Suppose first that $V_1 \oplus V_2$ is trivial. Then Lemma~\ref{diamonds} tells us that $V = V_0 \oplus V_3$ as an $\Out{3}$-module.
Suppose now that $V_0 \oplus V_3$ is trivial. Lemma~\ref{lem: dim 6 new case} tells us that $V = V_1 \oplus V_2$ as an $\Out{3}$-module. In both situations we can apply Lemma~\ref{lem: commuting delta}.

We are left with the most general case:
suppose that $\dim V_1 \oplus V_2 = 3$. We are going to assume that in fact $V_2 = \{0\}$, the other case being analogous. Applying Lemma~\ref{lem: gap in levels} gives us two $\Out 3$-representations $r: \Out 3 \to  V / (V_1 \oplus V_3)$ and $s: \Out 3 \to V / V_0$, where at least one of them occurs as a submodule of $V$.
Also, $r$ and $s$ factor through $\pi$ by Lemma~\ref{lem: commuting delta}. If any of these representations has dimension 0, then we are done. In what follows we shall suppose that the dimension of both $r$ and $s$ is non-zero, and thus that $V$ is reducible as an $\Out 3$-module.
We can choose a basis for $V$ so that the matrices in $\phi( \Out 3 )$ are all in a block-upper-triangular form, with diagonal blocks corresponding to representations $r$ and $s$.

Let $\Torelli_3 = \ker \pi$ be the Torelli subgroup. Our aim is to show that $\Torelli_3$ lies in the kernel of $\phi$.

Elements in $\Torelli_3$ map to matrices with identities on the diagonal, and all non-zero entries located in the block in the top-right corner. Hence $\Torelli_3$ maps to an abelian group isomorphic to $\K^m$, where $m \in \{5,8,9\}$ depends on the dimension of $r$.

Note that all products $\epsilon_i \epsilon_j$ lie in the kernel of $r$, and hence so do all elements \[\rho_{kj}^{2}=(\rho_{kj}^{\epsilon_i \epsilon_j} \rho_{kj}^{-1})^{-1},\]
where we took $k \neq i$. The work of Mennicke~\cite{mennicke1965} now shows that in fact $r$ factors through a finite group: when restricted to $\mathrm{SOut}(F_3) = \pi^{-1}(\SL_3(\Z))$, it factors through
\[\mathrm{SOut}(F_3) \to \SL_3(\Z) \to \SL_2(\Z_2).\] Let $A$ denote the kernel of this map.
Note that $\Torelli_3 < A$.

%

We have shown above that $r\vert_A$ is trivial,
and so $A$ maps to the the identity matrix in the block corresponding to $r$. Note that $\phi(A)$ acts by conjugation on the abelian group of matrices with identity blocks on the diagonal, and a trivial block in the bottom-left corner. As remarked above, this group is isomorphic to $\K^m$. Each row or column (depending on which diagonal block corresponds to $r$) in the top-right corner corresponds to an $A$-submodule, and so the group $\K^m$ splits as an $A$-module into
\[ \K^5, 2. \K^4 , 3. \K^3 , 4. \K^2  \textrm{ or } 5 . \K ,\]
depending on the dimension of $r$, where the multiplicative notation indicates the number of direct summands.

Let $T = \Torelli_3 / [\Torelli_3, \Torelli_3]$ denote the abelianisation of the Torelli group seen as an $\Out 3 / \Torelli_3 = \GL_3(\Z)$-module, where the action is the one induced by the conjugation action $\Out 3 \curvearrowright \Torelli_3$.
The structure of this module is known (see Kawazumi~\cite{kawazumi2005}) -- it is the second symmetric power of the dual of the standard $\GL_3(\Z)$-module. After tensoring $T$ with $\K$, we can apply Proposition~\ref{prop: irreducibility of K-module}, and conclude that $T \otimes_\Z \K$ is an irreducible $A$-module of dimension 6.
By Schur's Lemma, if we have an $A$-equivariant quotient of $T$, it is either isomorphic to $T$ or equal to $\{ 0 \}$.

Now consider the action of $\phi(A)$ on $\phi(\Torelli_3)\otimes_\Z \K$ by conjugation. It is at the same time an equivariant quotient of an irreducible 6-dimensional module and a submodule of \[ \K^5, 2. \K^4 , 3. \K^3 , 4. \K^2  \textrm{ or } 5 . \K .\]
This implies that
the image of $\Torelli_3$ under $\phi$ is trivial. This finishes the proof.
\end{proof}
\end{thm}

\section{Small graphs with transitive automorphism groups}
\label{secgraphs}

In this section we will establish some lemmata concerning graphs of rank at most $5$ with groups $W_3$ and $G_3$ acting on them.

\begin{dfn}[Admissible graphs]
Let $X$ be a connected graph with no vertices of valence 2, and suppose we have a group $G$ acting on it. We say that $X$ is \emph{$G$-admissible} \iff there is no $G$-invariant non-trivial (i.e. with at least one edge) forest in $X$. We also say that $X$ is \emph{admissible} \iff it is $\mathrm{Aut}(X)$-admissible.
\end{dfn}

The following theorem is due to Marc Culler~\cite{culler1984}, Dmitri Khramtsov~\cite{khramtsov1985} and Bruno Zimmermann~\cite{zimmermann1996} (each independently).

\begin{thm}[Culler~\cite{culler1984}; Khramtsov~\cite{khramtsov1985}; Zimmermann~\cite{zimmermann1996}]
\label{theorem}
Let $n \geqslant 2$. Suppose \[G \into \Out{n}\] is an embedding of a finite group $G$. Then there exists a finite $G$-admissible graph $X$ of rank $n$ (with a fixed isomorphism $\pi_1X \cong F_n$), so that the composition
\[G \to \mathrm{Aut}(X) \to \Out{n}\]
is the given embedding.
\end{thm}

The reason for requiring $n \geqslant 2$ in the above theorem is that there are no admissible graphs of rank 1 (i.e. with the fundamental group $\Z$). Of course any finite subgroup of $\Out{1} \cong \Z / 2 \Z$ can be realised as an action on a graph with one vertex and one edge.

Since we will be dealing with homology of finite graphs quite frequently in this section, let us observe the following.

\begin{lem}
Let $X$ be a finite, oriented graph. Recall that Definition~\ref{graphdef} gives us two maps $\iota, \tau : E(X) \to V(X)$. We have the following identification
\begin{displaymath}
H_1(X,\C) \cong \Big\{ f : E(X) \to \C \mid \forall v \in V(X) :  \sum_{\iota(e)=v} f(e) = \sum_{\tau(e)=v} f(e) \Big\}.
\end{displaymath}
\end{lem}

We will often refer to each such function $f$ as a \emph{choice of weights} of edges in $X$.

Before proceeding any further, we need to introduce a concept of \emph{collapsing maps} of graphs.

\begin{dfn}[Collapsing map]
\label{collapsingmap}
Let $\pi :X \to X'$ be a surjective morphism of graphs $X$ and $X'$. We say that $\pi$ is a \emph{collapsing map} \iff
the preimages of points in $X'$ under $\pi$ are path-connected.
\end{dfn}

\begin{rmk}
Let us observe three facts:
\begin{enumerate}
\item for a given graph $X$, giving a subset of $E(X)$ which will be collapsed specifies a collapsing map $\pi$;
\item any collapsing map $\pi : X \to X'$ induces a surjective map on homology by pushing forward weights of edges which it does not collapse. We will often abuse notation and refer to this induced map as $\pi$;
\item if $\pi : X \to X'$ is a $G$-equivariant collapsing map (where $G$ is a group acting on $X$ and $X'$), then the induced map on homology is also $G$-equivariant.
\end{enumerate}
\end{rmk}

It is now time to focus on the main area of our interest here.

\begin{dfn}
Let $G$ be a group acting on a graph $X$, and let $e$ be an edge of $X$. We define $X_e$ to be the graph obtained from $X$ by collapsing all edges \emph{not contained} in the $G$-orbit of $e$.
\end{dfn}

Note that the action of $G$ on such an $X_e$ is edge-transitive.

\begin{lem}
\label{eqn lemma}
Suppose $X$ is a $G$-admissible graph of rank at most 5, where $G$ is a group, and $e$ is any edge of $X$. Then $X_e$ has no vertices of valence 1 or greater than 10 and satisfies
\begin{equation*}
\label{*}
8 \geqslant 2 v_2 + v_3 + 2 v_4 + 3 v_5 + 4 v_6 +5 v_7 + 6 v_8 + 7 v_9 + 8 v_{10} \tag{$\ast$}
\end{equation*}
where $v_i$ is the number of vertices of valence $i$ in $X_e$.
\begin{proof}
First note that there are no vertices of valence 1 in $X_e$, since they could only occur if there were separating edges in $X$. But $X$ is admissible, and so there are no such edges.

A simple Euler characteristic count yields
\[ 2(\mathrm{rank}(X_e)-1) \geqslant \sum_{i=3}^\infty (i-2)v_i \]
and hence in particular $v_i = 0$ for all $i > 10$, as $X_e$ has rank at most 5.

Since $X$ is admissible, each vertex of $X_e$ of valence two comes from collapsing a subgraph of $X$ which is not a tree, hence
\[ \mathrm{rank}(X_e) \leqslant 5 -v_2 \]
and the result follows.
\end{proof}
\end{lem}

We will now consider graphs satisfying \eqref{*} with a transitive action of {${W_3 \cong G_3}$} yielding particular representations on the $\C$-homology of the graph.

\begin{prop}
\label{graph prop}
Let $X$ be a graph of rank 5 on which $G \in \{ W_3,  G_3 \}$ acts so that the representation of $G$ on $V = H_1(X,\C)$ induced by the action decomposes as $V = V_0 \oplus V_1 \oplus V_2 \oplus V_3$, where
\begin{itemize}
\item if $G = W_3$ then the decomposition is the one described in Definition~\ref{repofW_n}, and, as $S_3$-modules, $V_1$ is a permutation, $V_2$ signed permutation, $V_0$ a sum of trivial and $V_3$ a sum of determinant representations;
\item if $G = G_3$ then $\Delta$ acts as identity on $V_0 \oplus V_2$ and as minus the identity on $V_1 \oplus V_3$, and as $S_4$-modules, $V_1$ is a standard, $V_2$ signed standard, $V_0$ a sum of trivial and $V_3$ a sum of determinant representations.
\end{itemize}
Then, there is a subgraph $Y \leqslant X$ isomorphic to a 3-rose, on which $G$ acts in such a way that, as an $S_3$-module (where $S_3<W_3 \cap G_3$),  $H_1(Y,\C)$ contains the standard representation.
\begin{proof}
Let $v \in V$ be a vector belonging to a standard representation of $S_3<G$. It is represented by a choice of weights on edges of $X$. Let $e$ be an edge with a non-zero weight. Then the image of $v$ in $H_1(X_e,\C)$ is non-trivial, and hence Schur's Lemma informs us that $H_1(X_e,\C)$ contains a standard $S_3$-module.

Let $Z=X_e$. Lemma~\ref{eqn lemma} tells us that $Z$ satisfies \eqref{*}. Also, since $G$ acts transitively on edges of $Z$, there are at most two vertex-orbits of this action, and hence in particular at most two values $v_i$ can be non-zero. Let us list all possible values of $v_i$, noting that $i v_i = j v_j$ if there are vertices of valence $i$ and $j$ in $Z$, and that $v_i$ must be even if $i$ is odd and there are only vertices of valence $i$ in $Z$. All possible cases are summarised in Figure~\ref{fig: case table}.

\begin{figure}
\caption{Case table}
\label{fig: case table}
\begin{tabular}{c|ccccccccc|c|c}
  Case number & $v_2$ & $v_3$ & $v_4$ & $v_5$ & $v_6$ & $v_7$ & $v_8$ & $v_9$ & $v_{10}$ & edges & rank\\ \hline
  (1) & 4 &  &  &  &  &  &  &  & & 4 & 1 \\
  (2) & 3 & 2 &  &  &  &  &  &  & & 6 & 2 \\
  (3) & 3 &  &  &  &  &  &  &  & & 3 & 1 \\
  (4) & 2 &  & 1 &  &  &  &  &  & & 4 &2 \\
  (5) & 2 &  &  &  &  &  &  &  & &2&1 \\
  (6) & 1 &  &  &  &  &  &  &  & &1&1 \\
  (7) &  & 8 &  &  &  &  &  &  & &12&5 \\
  (8) &  & 6 &  &  &  &  &  &  & &9&4 \\
  (9) &  & 4 &  &  &  &  &  &  & &6&3 \\
  (10) &  & 2 &  &  & 1 &  &  &  & &6&4 \\
  (11) &  & 2  &  &  &  &  &  &  & &3&2 \\
  (12) &  &  & 4 &  &  &  &  &  & &8&5 \\
  (13) &  &  & 3 &  &  &  &  &  & &6&4 \\
  (14) &  &  & 2 &  &  &  &  &  & &4&3 \\
  (15) &  &  & 1 &  &  &  &  &  & &2&2 \\
  (16) &  &  &  & 2 &  &  &  &  & &5&4 \\
  (17) &  &  &  &  & 2 &  &  &  & &6&5 \\
  (18) &  &  &  &  & 1 &  &  &  & &3&3 \\
  (19) &  &  &  &  &  &  & 1 &  & &4&4 \\
  (20) &  &  &  &  &  &  &  &  & 1 &5&5 \\
\end{tabular}
\end{figure}

Now, in order to have a standard \rep of $S_3$, we need at least 3 edges in $Z$, and the rank of $Z$ has to be at least 2. We can therefore immediately rule out cases $(1), (3), (5),(6)$ and  $(15)$. Also, since the action of $G$ on the edges of $Z$ is transitive, their number has to divide $\vert G \vert = 48$. Hence we can additionally rule out cases $(8), (16)$ and $(20)$. We are left with the cases listed in Figure~\ref{fig: reduced case table}.

\begin{figure}
\caption{Reduced case table}
\label{fig: reduced case table}
\begin{tabular}{c|ccccccccc|c|c}
  Case number & $v_2$ & $v_3$ & $v_4$ & $v_5$ & $v_6$ & $v_7$ & $v_8$ & $v_9$ & $v_{10}$ & edges & rank\\ \hline
  (2) & 3 & 2 &  &  &  &  &  &  & & 6 & 2 \\
  (4) & 2 &  & 1 &  &  &  &  &  & & 4 &2 \\
  (7) &  & 8 &  &  &  &  &  &  & &12&5 \\
  (9) &  & 4 &  &  &  &  &  &  & &6&3 \\
  (10) &  & 2 &  &  & 1 &  &  &  & &6&4 \\
  (11) &  & 2  &  &  &  &  &  &  & &3&2 \\
  (12) &  &  & 4 &  &  &  &  &  & &8&5 \\
  (13) &  &  & 3 &  &  &  &  &  & &6&4 \\
  (14) &  &  & 2 &  &  &  &  &  & &4&3 \\
  (17) &  &  &  &  & 2 &  &  &  & &6&5 \\
  (18) &  &  &  &  & 1 &  &  &  & &3&3 \\
  (19) &  &  &  &  &  &  & 1 &  & &4&4 \\
\end{tabular}
\end{figure}

We will need to deal with these cases one by one:

\noindent \textbf{Case (2):} Here we have three vertices of valence two, on which $S_3$ has to act transitively. Each of these comes from collapsing a graph of non-zero rank in $X$, hence the sum of homologies of these graphs contains another standard \rep of $S_3$. This contradicts our assumptions.

\noindent \textbf{Case (4):} Here $Z$ is a subdivided 2-rose, so we cannot get a standard \rep of $S_3$ on the homology of this graph.

\noindent \textbf{Case (7):} There are four graphs with an edge-transitive group action with at most 8 vertices each of valence 3, namely a 3-cage, the 1-skeleton of a tetrahedron, the complete bipartite graph $K(3,3)$, and the 1-skeleton of a cube.
Clearly, only the last one has the required number of vertices.
An $S_3$-action yielding a standard \rep has to be the one given by fixing two vertices and permuting 3 edges incident at one of them in a natural way. This however yields two copies of the standard \rep when acting on homology, which contradicts our assumptions.

\noindent \textbf{Case (9):} In this case we are dealing with an edge-transitive $G$-action on 1-skeleton of a tetrahedron. Such an action has to also be vertex-transitive, but $S_3$ cannot act transitively on 4 points.

%
%
\noindent \textbf{Case (10):} In this case $Z$ is obtained by taking a wedge of two 3-cage graphs, $C_1$ and $C_2$, say. Since the action of $G$ is edge-transitive, for any edge $e$ in $Z$ we have
\[ \vert \mathrm{Stab}_G(e) \vert = 8\]
and hence each 3-cycle acts freely. Since a 3-cycle cannot swap $C_1$ and $C_2$, it must act in the natural way on edges of both. Hence there have to be two copies of the standard \rep of $S_3$ in $H_1(\Z,\C)$, which is not the case.

%
%
%
%
%
\noindent \textbf{Case (11):} In this case $Z$ is a 3-cage. If $G = G_3$, then $S_4 < G_3$ cannot act on $Z$ yielding the desired standard or signed standard representation. Suppose now that $G = W_3$. If $\epsilon_1$ preserves exactly one edge, then so do $\epsilon_2$ and $\epsilon_3$; these edges are distinct, as otherwise we would have some $\epsilon_i$ and $\epsilon_j$ acting in the same way where $i\neq j$, and so $H_1(Z,\C) \leqslant V_0 \oplus V_3$, where $S_3$ cannot have a standard representation. Since the edges are distinct, $\epsilon_1$ and $\epsilon_2$ cannot commute. This shows that each $\epsilon_i$ preserves all edges of $Z$, and hence $H_1(Z,\C) \leqslant V_0 \oplus V_3$, which is a contradiction.

\noindent \textbf{Case (12):} The graph $Z$ is a bipartite graph with 4 vertices and exactly zero or two edges connecting each pair of vertices. Hence $Z$
admits a $G$-equivariant quotient map to a square (i.e. a single cycle made of 4 edges) which is a $2$-to-$1$ map on edges. Each 3-cycle acts trivially on the square; moreover it cannot act non-trivially on the preimages of edges of the square. We conclude that each 3-cycle acts trivially, which is a contradiction.

%
%
\noindent \textbf{Case (13):} We easily check that the graph $Z$ consists of three vertices, each of which has exactly two edges connecting it to each of the other two. Since the 3-cycle in $S_3$ acts non-trivially, it has to act transitively on vertices, and so either each vertex in $Z$ comes from collapsing a subgraph which was not a tree in $X$, or none of them does. Neither of these two cases is possible, since the rank of $X$ is 5.

\noindent \textbf{Case (14):} In this case $Z$ is a 4-cage. Each edge in $Z$ has a corresponding edge in $X$, and the fact that $X$ is $G$-admissible implies that these edges do not form a forest. Hence they can form either  a single simple loop, a pair of simple loops, or a 4-cage in $X$. The first two cases are impossible, since they would yield a trivial action of the 3-cycle of $S_3$ on $Z$. Hence $X$ contains $Z$ as a subgraph.

Our assumption on the \reps of $G$ tells us that either $\Delta$ or each transposition in $S_3$ has to flip $Z$, and so $X$ is a 4-cage with a loop of length one attached to each vertex. Let $a$ and $b$ be two vectors in $H_1(X, \C)$, each given by putting a weight 1 on exactly one of the loops.

If $\Delta$ flips the graph, then $a+b$ and $a-b$ span two one-dimensional eigenspaces of $\Delta$, one with eigenvalue +1, and one with eigenvalue $-1$. Hence transpositions in $S_3$ have to map one of these vectors to itself, and the other to minus itself; this is only possible if they flip the graph, which contradicts our assumptions.

A similar argument works if the transpositions in $S_3$ flip $X$.

\noindent \textbf{Case (17):}
In this case $Z=X$ is a 6-cage. As before we have
\[\vert \mathrm{Stab}_G(e) \vert = 8\]
for any edge $e$ in $Z$. Hence each 3-cycle in $S_3$ acts freely and so we have two copies of the standard $S_3$-representation, which is a contradiction.


\noindent \textbf{Case (18):} In this case $Z$ is a 3-rose. If $Z$ is actually a subgraph of $X$, then we are done. Suppose it is not.

As $Z$  only has one vertex, there is a connected subgraph $X'$ of $X$ that we collapsed when constructing $Z$. Since $X'$ is of rank 2, after erasing vertices of valence 2 (in $X'$), we are left with two cases: a 2-rose or a 3-cage. Since $Z$ is not a subgraph of $X$, and the preimages of edges of $Z$ in $X$ cannot form a forest, they either form a simple loop (of length three), or a disjoint union of three loops (each of length one). In any event, we have three vertices on which the 3-cycle in $S_3$ acts transitively. Hence $X'$ has to be a 3-cage, with the three vertices lying on the three edges of the cage. But then we get two standard \reps of $S_3$ inside the homology of $X$, which is a contradiction.

\noindent \textbf{Case (19):} In this case $Z$ is a 4-rose. The 3-cycle in $S_3$ acts by permuting three petals, and fixing one; let us call this fixed edge $f$. We easily check that $f$ is preserved setwise by $S_3$, and hence also by $\Delta$, since $\Delta$ commutes with $S_3$.

If $G = W_3$ then the one-dimensional subspace in $H_1(Z ,\C)$ spanned by a vector corresponding to $f$ is contained either in $V_0$ or in $V_3$, and hence $f$ has to be preserved by all elements in $G$. This contradicts transitivity of the action of $G$ on $Z$.

Suppose $G = G_3$. Note that there is only one way (up to isomorphism) in which $S_4$ can act on a set of four elements transitively. Therefore, as $\Delta$ commutes with $S_4$, $\Delta$ acts as plus or minus the identity on $H_1(Z,\C)$. Now $H_1(Z,\C)$ as an $S_4$-\rep is a sum of standard and trivial or signed standard and determinant representations. In particular, our hypothesis tells us that $\Delta$ cannot act as either plus or minus the identity.
This is a contradiction.
\end{proof}
\end{prop}

\begin{lem}
\label{trivial case}
Suppose $\phi : \Out{3} \to \Out{5}$ is a homomorphisms. Let \[\psi : \Out{5} \to \GL(V) \cong \GL_5(\C)\] be the natural map, where $V=H_1(F_5,\C)$. Note that $\psi \circ \phi$ gives a representation of $\Out{3}$ on $V$. If, as a $W_n$-module, $V$ splits as $V_0 \oplus V_3$ (with the notation of Definition~\ref{repofW_n}) then the image of $\phi$ is finite.
\begin{proof}
The fact that $V=V_0 \oplus V_3$ as a $W_n$-module implies that $\psi \circ \phi (\epsilon_i \Delta) =1$ for each $i$. Now a result of Baumslag--Taylor~\cite{baumslagtaylor1968} tells us that the kernel of $\psi$ is torsion-free, and so $\phi(\epsilon_i \Delta)=1$. But this means that we have the following commutative diagram:
\[ \xymatrix{
\Out 3 \ar[d] \ar[r]^{\phi} & \Out 5 \\
\Out{3}/\langle \! \langle \{ \rho_{ij}=\lambda_{ij} : i \neq j\}\rangle \! \rangle \ar[ur] \ar[r]^{\ \ \ \ \ \ \ \ \ \ \  \ \ \ \ \cong} & \GL_3(\Z) \ar[u]_{\phi'}
} \]
This allows us to use a result of Bridson and Farb~\cite{bridsonfarb2001}, who have shown that such a $\phi'$ necessarily has finite image. Therefore the image of $\phi$ is finite.
\end{proof}
\end{lem}

\begin{thm}
\label{result C}
\label{main result}
Suppose $\phi : \Out{3} \to \Out{5}$ is a homomorphism. Then the image of $\phi$ is finite.
\begin{proof}
Consider the natural map $\psi : \Out{5} \to \GL_5(\C) \cong \GL(V)$ as above. Again as above, the composition $\eta = \phi \circ \psi$ gives us a 5-dimensional complex linear \rep $\eta : \Out{3} \to \GL_5(\C)$.

Suppose first that, with the notation of Definition~\ref{repofW_n}, $V$ satisfies
\[V = V_0 \oplus V_3 \]
 Then Lemma~\ref{trivial case} yields the result.

Now suppose that $V = V_0 \oplus V_1 \oplus V_2 \oplus V_3$ where $V_1 \oplus V_2 \neq \{ 0 \}$. We apply Proposition~\ref{reps prop} to $\eta$.

We will now apply Theorem~\ref{theorem} to two finite subgroups of $\Out{5}$, namely $\phi(G_3)$ and $\phi(W_3)$, to obtain two graphs $X$ and $Y$ respectively, on which the groups $G_3$ and $W_3$ act. Note that $H_1(X, \C) \cong H_1(Y, \C) \cong V$, and the \reps of $G_3$ and $W_3$ induced by the actions of the groups on homology of the respective graphs are isomorphic to the ones given by restricting $\eta$. Hence the conclusions of Proposition~\ref{reps prop} apply to these representations, and so we can apply Proposition~\ref{graph prop} to the actions $G_3 \curvearrowright X$ and $W_3 \curvearrowright Y$.

We conclude that both $X$ and $Y$ have a subgraph, preserved by the action of the respective group, isomorphic to a 3-rose. We also know that we can label the petals as $e_1, e_2, e_3$, so that $S_3$ acts on this rose by permuting petals in the natural way, with the transpositions potentially also flipping all petals.

Knowing that $V_1 \oplus V_2 \neq \{ 0 \}$ implies that in the $W_3$ case, either each $\epsilon_i$ flips $e_i$ and leaves the other petals fixed, or each $\epsilon_i$ fixes $e_i$ and flips the other petals. In the $G_3$ case, we see that there is only one way in which $S_4$ can act on the 3-rose inducing a standard or a signed standard representation. Each $\sigma_{i4}$ has to interchange the two petals with labels different than $e_i$ and preserve the third one; additionally, it either flips $e_i$ and keeps some orientation of the other two fixed, or it flips the other two and fixes $e_i$. These two cases depend on the action of $\sigma_{ij}$ for $i,j \leqslant 3$.

In any case, we have
\[ \phi(\sigma_{14}) = \phi(\sigma_{23} \epsilon_2 \epsilon_3) \]
and so
\[ \phi(\lambda_{21}) = \phi(\lambda_{21}^{\sigma_{14} \epsilon_1}) = \phi(\lambda_{21})^{\phi(\Delta \sigma_{23})} = \phi(\rho_{31}) \]
Therefore
\[ 1 = \phi([\rho_{23}^{-1}, \lambda_{21}^{-1}]) = \phi([\rho_{23}^{-1}, \rho_{31}^{-1}]) = \phi(\rho_{21})^{-1} \]
It follows that all Nielsen moves (which generate an index 2 subgroup of $\Out 3$) lie in the kernel od $\phi$, and so the image of $\phi$ is of size at most 2, determined by $\phi(\epsilon_1)$.
\end{proof}
\end{thm}

%
\bibliographystyle{abbrv}
\bibliography{bibliography}

\begin{thebibliography}{10}

\bibitem{aramayonaetal2009}
J.~Aramayona, C.~J. Leininger, and J.~Souto.
\newblock Injections of mapping class groups.
\newblock {\em Geom. Topol.}, 13(5):2523--2541, 2009.

\bibitem{baumslagtaylor1968}
G.~Baumslag and T.~Taylor.
\newblock The centre of groups with one defining relator.
\newblock {\em Math. Ann.}, 175:315--319, 1968.

\bibitem{bogopol'skiipuga2002}
O.~Bogopol'ski{\u\i} and D.~Puga.
\newblock On embeddings of $\mathrm{Out}({F}_n)$, the outer automorphism group
  of the free group of rank $n$, into $\mathrm{Out}({F}_m)$ for $m > n$.
\newblock {\em Algebra and Logic}, 41(2):69--73, 2002.

\bibitem{bridsonfarb2001}
M.~R. Bridson and B.~Farb.
\newblock A remark about actions of lattices on free groups.
\newblock {\em Topology Appl.}, 110(1):21--24, 2001.

\bibitem{bridsonvogtmann2003}
M.~R. Bridson and K.~Vogtmann.
\newblock Homomorphisms from automorphism groups of free groups.
\newblock {\em Bull. London Math. Soc.}, 35(6):785--792, 2003.

\bibitem{bridsonvogtmann2011}
M.~R. Bridson and K.~Vogtmann.
\newblock Abelian covers of graphs and maps between outer automorphism groups
  of free groups.
\newblock {\em Mathematische Annalen}, pages 1--34, 2011.

\bibitem{culler1984}
M.~Culler.
\newblock Finite groups of outer automorphisms of a free group.
\newblock In {\em Contributions to group theory}, volume~33 of {\em Contemp.
  Math.}, pages 197--207. Amer. Math. Soc., Providence, RI, 1984.

\bibitem{gersten1984}
S.~M. Gersten.
\newblock A presentation for the special automorphism group of a free group.
\newblock {\em J. Pure Appl. Algebra}, 33(3):269--279, 1984.

\bibitem{grunewaldlubotzky2006}
F.~Grunewald and A.~Lubotzky.
\newblock Linear representations of the automorphism group of a free group.
\newblock {\em Geometric and Functional Analysis}, 18(5):1564--1608, June 2006.

\bibitem{kawazumi2005}
N.~Kawazumi.
\newblock Cohomological aspects of magnus expansions.
\newblock {\em {arXiv}:0505497}.

\bibitem{khramtsov1985}
D.~G. Khramtsov.
\newblock Finite groups of automorphisms of free groups.
\newblock {\em Mat. Zametki}, 38(3):386--392, 476, 1985.

\bibitem{khramtsov1990}
D.~G. Khramtsov.
\newblock Outer automorphisms of free groups.
\newblock In {\em Group-theoretic investigations ({R}ussian)}, pages 95--127.
  Akad. Nauk SSSR Ural. Otdel., Sverdlovsk, 1990.

\bibitem{kielak2013}
D.~Kielak.
\newblock Outer automorphism groups of free groups: linear and free
  representations.
\newblock {\em J. London Math. Soc.}, 87(3):917--942, 2013.

\bibitem{mennicke1965}
J.~L. Mennicke.
\newblock Finite factor groups of the unimodular group.
\newblock {\em Ann. of Math. (2)}, 81:31--37, 1965.

\bibitem{potapchikrapinchuk2000}
A.~Potapchik and A.~Rapinchuk.
\newblock Low-dimensional linear representations of {${\rm Aut}\,F_n,\ n\geq
  3$}.
\newblock {\em Trans. Amer. Math. Soc.}, 352(3):1437--1451, 2000.

\bibitem{TurchinWillwacher2015}
V.~Turchin and T.~Willwacher.
\newblock Hochschild--{P}irashvili homology on suspensions and representations
  of $\mathrm{Out}({F}_n)$.
\newblock {\em In preparation}.

\bibitem{zimmermann1996}
B.~Zimmermann.
\newblock Finite groups of outer automorphisms of free groups.
\newblock {\em Glasgow Math. J.}, 38(3):275--282, 1996.

\end{thebibliography}
\end{document}